\documentclass[num-refs]{wiley-article}

\usepackage{graphicx}
\usepackage[space]{grffile}
\usepackage{latexsym}
\usepackage{textcomp}
\usepackage{longtable}
\usepackage{tabulary}
\usepackage{booktabs,array,multirow}
\usepackage{amsfonts,amsmath,amssymb}
\usepackage{natbib}
\usepackage{url}
\usepackage{hyperref}
\hypersetup{colorlinks=false,pdfborder={0 0 0}}
\usepackage{etoolbox}
\usepackage{pdfpages}
\usepackage{subcaption}
\newif\iflatexml\latexmlfalse

\AtBeginDocument{\DeclareGraphicsExtensions{.pdf,.PDF,.eps,.EPS,.png,.PNG,.tif,.TIF,.jpg,.JPG,.jpeg,.JPEG}}

\usepackage[utf8]{inputenc}
\usepackage[english]{babel}

\usepackage{siunitx}
\usepackage{microtype}  
\usepackage{titlesec}
\titleformat{\section}
  {\normalfont\fontsize{10}{12}\bfseries\raggedright}
  {\thesection}{1em}{\textmd{|}\quad\MakeUppercase{#1}}

\iflatexml

\else

\paperfield{Systems Engineering}
\abbrevs{DSM, Design Structure Matrix; MDM, Multi-Domain Matrix; HR, Human Resources; HRM, Human Resource Management; LLM, Large Language Model; SME, Small and Medium Enterprise; UI/UX, User Interface/User Experience.}
\corraddress{Jaewoo Kim, Department of Aerospace Engineering, Korea Advanced Institute of Science and Technology, 291 Daehak-ro, Daejeon 34141, Republic of Korea.}
\corremail{jw.kim@kaist.ac.kr}
\fundinginfo{}
\fi

\papertype{Regular Article}

\title{Multi-Domain Matrix Framework for Human Resource Decision Support}

\author[1,2]{Daehan Lee}
\author[1]{Jaewoo Kim}
\author[1]{Jaemyung Ahn}

\affil[1]{Department of Aerospace Engineering, Korea Advanced Institute of Science and Technology, Daejeon, Korea.}
\affil[2]{Planby Technologies, Inc., Seoul, Korea.}

\runningauthor{Daehan Lee}

\begin{document}

\maketitle
\selectlanguage{english}
\begin{abstract}
This paper presents an actionable human resource (HR) decision-support framework for small firms and startups based on a multi-domain matrix (MDM). The framework addresses three key challenges faced by small organizations: complex interdependencies among organizational components; the lack of systematic analytical tools for HR decision-making; and the need for rapid responses in fast-changing organizational environments. The proposed framework formulates startup human resource management as a multi-domain structural modeling problem, where members, skills, and projects are interconnected domains within an integrated MDM. Based on this representation, the framework provides qualitative analysis guidelines and quantitative metrics for diagnosing an organization's HR state and supporting personnel decisions on workload redistribution, hiring, and capability development. A case study of MDM-based HR decisions for an early-stage technology startup is conducted to demonstrate the framework's practical applicability. The application shows that the framework can identify workload imbalances, reveal a key member with an unsustainable workload, and inform a subsequent hiring decision. The framework can be further applied after the hiring of a new member to track changes in the organization's multi-domain structure and support continuous HR diagnosis.

\end{abstract}%

\twocolumn

\section{Introduction}
\label{sec:introduction}

Human resources (HR) refers to the people who constitute an organization's workforce, together with their collective knowledge, skills, abilities, and experiences. From a resource-based view, such intangible and firm-specific resources can contribute to sustained competitive advantage~\cite{barney1991firm,wright1994human}. Reflecting this importance, human resource management (HRM) encompasses the strategic management of workforce-related activities, including recruitment, training, performance management, and retention, with the aim of aligning human capital with organizational objectives~\cite{wright1992theoretical}.

HRM is one of the most critical challenges facing early-stage startups. Unlike large firms with dedicated HR departments and established processes, startups must make decisions under conditions of limited resources, incomplete information, and high uncertainty~\cite{cardon2004managing}. Although this study focuses on early-stage startups, it draws on the broader literature on small and medium enterprises (SMEs) and small entrepreneurial firms because these organizations share several relevant characteristics, including resource constraints, limited functional specialization, informal HR practices, and high dependence on a small number of individuals~\cite{heneman2000human,cardon2004managing,vanlancker2022hrm}. These difficulties are further intensified by the structural and operational characteristics of startups, where organizational actors, capabilities, responsibilities, activities, and resources are tightly interdependent, yet such interdependencies are often not directly visible to decision-makers.

Under these conditions, wrong HR decisions can have especially severe consequences. Because startups operate with limited slack, even a single poor HR decision can disproportionately affect their performance and survival. More broadly, slow or poorly informed HR decisions can reduce adaptability, lose opportunities, and, in some cases, jeopardize organizational survival.

Motivated by this gap, this study proposes a framework for organizational analysis and HR decision support for early-stage startups. The framework uses the multi-domain matrix (MDM), a systems engineering tool for representing interdependencies across multiple domains, as its core analytical tool. Rather than treating HR decisions as isolated personnel-level judgments, the proposed framework models the organization as a set of interrelated domains, such as organizational actors, capabilities, responsibilities, activities, and resources. These domains and their cross-domain relationships are represented within an integrated MDM structure, enabling decision-makers to examine how local personnel decisions may affect the broader organizational system. In the case study presented in this paper, this general multi-domain logic is instantiated through three HR-relevant domains: members, skills, and projects.

The contribution of this study is threefold. First, we formulate the HRM as a multi-domain interdependency problem and develop an MDM-based representation to capture the relationships among HR-relevant organizational domains. Second, we derive interpretable quantitative indicators from this representation, including the Workload Index and Value Index, to assess workload concentration, individual contribution, burnout risk, underutilization, hiring priorities, and capability development needs. Third, we demonstrate how the proposed framework can support actionable HR decisions in an early-stage startup context through a case study. Thus, the proposed MDM-based framework provides startups with a systematic yet practical approach for capturing complex organizational interactions while remaining readily adaptable to different domain configurations and rapidly changing organizational environments.

The remainder of this paper is organized as follows. Section~\ref{sec:literature} reviews relevant literature on HR challenges in startups, DSM/MDM methodologies, and the research gap addressed in this study. Section~\ref{sec:framework} presents the proposed framework, including the general MDM-based modeling procedure, domain instantiation, metric formulations, and the quadrant-based analytical approach. Section~\ref{sec:case-study} applies the framework to an early-stage startup and discusses the resulting organizational insights and HR implications. Section~\ref{sec:conclusions} concludes the paper by summarizing the main contributions, discussing limitations, and outlining directions for future research.

\section{Literature Review}
\label{sec:literature}

\subsection{HR Challenges in Startups}
\label{subsec:hr-startups}

Human resource management (HRM) in startups presents distinctive challenges that differ substantially from those in established organizations. Startups and small entrepreneurial firms typically operate with limited resources, low functional specialization, and high dependence on a small number of individuals~\cite{heneman2000human,cardon2004managing,vanlancker2022hrm}. In such contexts, employees rarely perform only a single well-defined role. Instead, individual members often contribute across multiple responsibilities, capability areas, and informal organizational functions. Their organizational value may therefore arise not only from assigned project tasks, but also from relational contributions such as knowledge sharing, mentoring, coordination, and communication facilitation. From a network-theoretic perspective, such contributions can be understood as arising from actors' positions and patterns of ties within the organization~\cite{borgatti2011network}.

These internal complexities are further intensified by the external conditions under which startups operate. Early-stage firms must adapt to changing market demands, technological uncertainty, investor expectations, and rapid scaling requirements. Prior research has shown that HR problems in small and medium-sized firms vary across organizational growth stages, indicating that people-management issues evolve as firms expand rather than remaining static~\cite{rutherford2003human}. Studies on entrepreneurial scaling also emphasize that growth creates major challenges in internal organizing, including changes in roles, routines, coordination mechanisms, and leadership responsibilities~\cite{desantola2017scaling,vanlancker2023preparing}. As a result, HR decisions in startups are not merely administrative decisions; they are time-sensitive organizational design decisions that affect whether the firm can sustain growth.

Despite this importance, startups often lack formalized HRM systems and analytical decision-support mechanisms. HR practices in small and entrepreneurial firms are frequently informal and founder-driven rather than supported by formal HR infrastructures~\cite{cardon2004managing,baron2002organizational,davila2005exploratory}. Moreover, SME owner-managers often handle HRM issues reactively, particularly when such issues become acute managerial concerns~\cite{tocher2009perceived}. Although these informal practices may provide flexibility in the early stage, they also make it difficult to diagnose workload imbalance, capability gaps, role ambiguity, burnout risk, and key-person dependency in a repeatable manner. Consequently, decisions regarding hiring, task allocation, workload redistribution, and capability development often rely heavily on managerial intuition. This reliance on intuition becomes problematic when organizational interdependencies are complex, because decision-makers may fail to identify hidden bottlenecks, overloaded personnel, or underutilized capabilities.

The consequences of poor HRM are particularly severe in startups. Because startups have limited slack resources, they are less able to absorb hiring mistakes, turnover, or the loss of key personnel. Studies on young high-technology firms show that changes in organizational employment models can increase employee turnover and negatively affect organizational performance subsequently~\cite{baron2001labor,baron2002organizational}. Moreover, non-founder human capital has been shown to affect the long-run growth and survival of high-technology ventures, indicating that the availability and allocation of complementary capabilities can be central to venture survival~\cite{siepel2017nonfounder}. From this perspective, HRM failure in startups is not merely an internal management inefficiency; it can directly threaten business continuity. In particular, when critical knowledge and responsibilities are concentrated in a few individuals, excessive job demands may increase burnout risk, thereby amplifying organizational fragility in resource-constrained startups~\cite{demerouti2001job}.

These challenges suggest a need for a practical analytics tool in startup HRM. However, quantitative and structured tools tailored to resource-constrained small organizations remain limited. The HR analytics studies have multiple limitations to HR analytics adoption, including issues related to data and models, software and technology, people, and management~\cite{fernandez2021tackling}. Recent reviews have also clarified determinants of effective HR analytics implementation and emphasized the growing importance of people analytics~\cite{wang2024determinants,yoon2024people}. Nevertheless, these studies primarily focus on the adoption and implementation conditions of HR analytics, rather than on how actionable insights can be systematically derived from small-scale organizational data. Moreover, these studies are generally not designed specifically for resource-constrained startup contexts, where HR data are sparse, roles are fluid, and formal HR infrastructure is limited. Taken together, the literature survey points to a clear gap: the need for a practical and systematic analytical framework that can support HR decision-making in startups under conditions of internal complexity, rapid adaptation, and limited analytical infrastructure.

\subsection{DSM/MDM-Based Organizational Analysis}
\label{subsec:dsm-org-analysis}

Design structure matrix (DSM) is a matrix-based representation for modeling dependencies among elements within a complex system, originally introduced by Steward as a method for managing the design of complex systems~\cite{steward1981design}. In a DSM, the elements of a system are arranged in the same order along the rows and columns of a square matrix, and the presence or strength of relationships between elements is recorded at each intersection. Its principal advantage is that it enables concise visualization and quantitative analysis of complex interdependencies. Owing to these properties, DSM has been widely recognized as a useful tool for modeling and analyzing complex systems, including product, process, project, and organizational structures~\cite{browning2001applying,browning2016dsm}.

\begin{figure}[h]
    \centering
        \includegraphics[page=13, width=0.8\columnwidth]{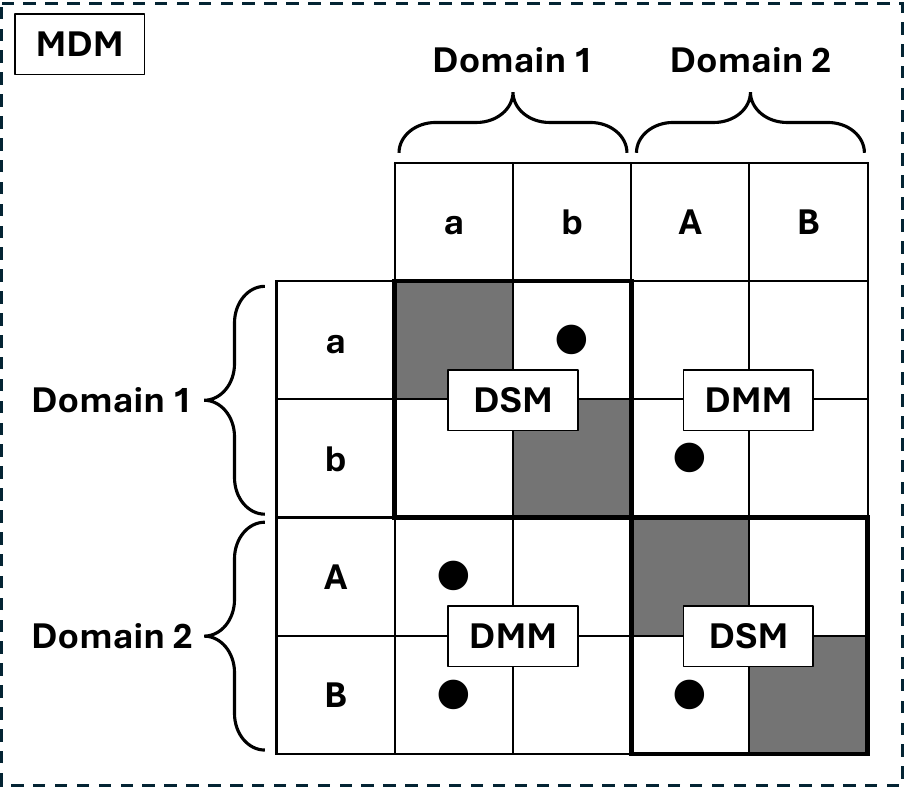}
    \caption{Structure of an MDM: diagonal same-domain blocks correspond to DSMs, whereas off-diagonal cross-domain blocks correspond to DMMs. Dots indicate the existence of relations between pairs of elements.}
    \label{fig:dsm-mdm}
\end{figure}

MDM extends DSM by representing interactions across multiple domains rather than within a single domain alone~\cite{danilovic2007managing,bartolomei2012engineering}. In an MDM, same-domain blocks correspond to DSMs, whereas cross-domain blocks are referred to as domain mapping matrices (DMMs). Bartolomei et al.~\cite{bartolomei2012engineering} proposed the Engineering Systems Multiple-Domain Matrix as an organizing framework for modeling large-scale complex systems, emphasizing its usefulness for representing heterogeneous system elements and their interactions within a unified structure. Related MDM methodologies have also been formalized within the Structural Complexity Management framework, which provides systematic procedures for modeling, analyzing, and optimizing multi-domain systems~\cite{lindemann2009structural}. Compared with a traditional single-domain DSM, MDM enables a more comprehensive representation of systems by explicitly incorporating relationships among domains of different natures. This feature is particularly relevant for organizational analysis because prior work has shown that technical architectures, development processes, and organizational structures can be analyzed as interdependent network structures~\cite{sosa2004network,parraguez2016characterizing}. Figure~\ref{fig:dsm-mdm} illustrates an MDM with two domains: one containing components a and b, and the other containing components A and B. The resulting matrix consists of two DSM blocks and two DMM blocks.

Prior studies have applied matrix- and network-based perspectives to product-development and organizational analysis. Eppinger et al.~\cite{eppinger1994model} proposed a matrix-based method for organizing product development tasks by capturing task sequences and technical relationships. In the systems engineering literature, Collins et al.~\cite{collins2009evaluating} integrated DSM and network analysis to evaluate task interactions in product development processes, showing that network metrics can help identify important tasks, critical interactions, and execution constraints. Browning~\cite{browning1999designing} further emphasized the importance of explicitly considering organizational integration requirements when designing system development projects.

Related work has also examined engineering organizations and teams as socio-technical networks. Avnet~\cite{avnet2016network} applied a network-based approach to analyze team coordination and shared cognition in systems engineering design teams. Parraguez et al.~\cite{parraguez2016characterizing} characterized design process interfaces as organization networks composed of people and their interactions, demonstrating how structural and compositional aspects of such networks can provide insights for engineering systems management. At a more general theoretical level, Baldwin and Clark~\cite{baldwin2000design} provided a foundational modularity perspective for decomposing complex systems into relatively independent modules governed by design rules. Together, these studies suggest that structural modeling approaches are well-suited to uncovering hidden patterns that may influence organizational coordination, decision-making, and performance.

Despite these advances, the application of DSM- and MDM-based methods to HR decision support remains limited. Existing systems engineering applications of DSM, network analysis, and organization modeling have largely emphasized product development processes, design interfaces, coordination, communication, organizational integration, and design structure, rather than integrated support for startup HR decisions. In particular, to the best of the authors' knowledge, prior studies have not developed an MDM-based framework specifically for startup HR decision support that translates cross-domain organizational relationships into concrete personnel decisions, including hiring, workload redistribution, and capability development. Existing studies provide useful foundations for structural representation and organizational network analysis, but they do not offer a lightweight and configurable framework for deriving HR decision-support metrics from multi-domain organizational data in resource-constrained startup settings.

\subsection{Research Gap Analysis and Contribution}
\label{subsec:research-gap}

In summary, the literature reveals a specific but unresolved gap at the intersection of HRM, HR analytics, and structural system modeling---particularly for small businesses like startups. The startup HRM literature shows that personnel decisions in early-stage firms are shaped by high interdependence, environmental dynamism, and limited analytical infrastructure, yet it offers limited methodological guidance for translating small-scale organizational data into decision-relevant insights. The HR analytics literature has advanced understanding of adoption conditions, barriers, and organizational readiness~\cite{fernandez2021tackling,wang2024determinants,yoon2024people}, but remains comparatively weak in providing lightweight analytical frameworks that early-stage startups can use to support concrete personnel decisions. In parallel, DSM/MDM and systems engineering research offer mature methodologies for representing and analyzing interdependent systems, including product development processes, organizational integration, team coordination, and engineering design interfaces~\cite{browning2016dsm,collins2009evaluating,browning1999designing,avnet2016network,parraguez2016characterizing}. However, these studies have not directly addressed HR decision support in early-stage startup contexts. Adjacent workforce-risk studies provide useful constructs for assessing burnout and job-related strain, but they do not by themselves represent HR-related organizational data as an interconnected multi-domain system~\cite{maslach2008early,demerouti2001job}.

Accordingly, the literature lacks a practical framework that integrates these streams. The missing contribution is therefore not merely the application of HR analytics to startups, nor simply the use of MDM in an organizational context, but the development of a quantitative decision-support framework that connects multi-domain structural modeling with concrete HR actions. In this study, the framework is demonstrated using members, skills, and projects, but the underlying modeling logic is not restricted to this particular domain configuration.

This study addresses that gap by developing an MDM-based HR decision-support framework for early-stage startups. In doing so, it contributes methodologically by extending MDM from structural representation to HR decision support, and contributes practically by providing a systematic yet implementable analytical tool for small, resource-constrained organizations.

\section{MDM-Based HR Decision Support Framework}
\label{sec:framework}

\begin{figure}[h]
    \centering
        \includegraphics[page=3, width=0.99\columnwidth]{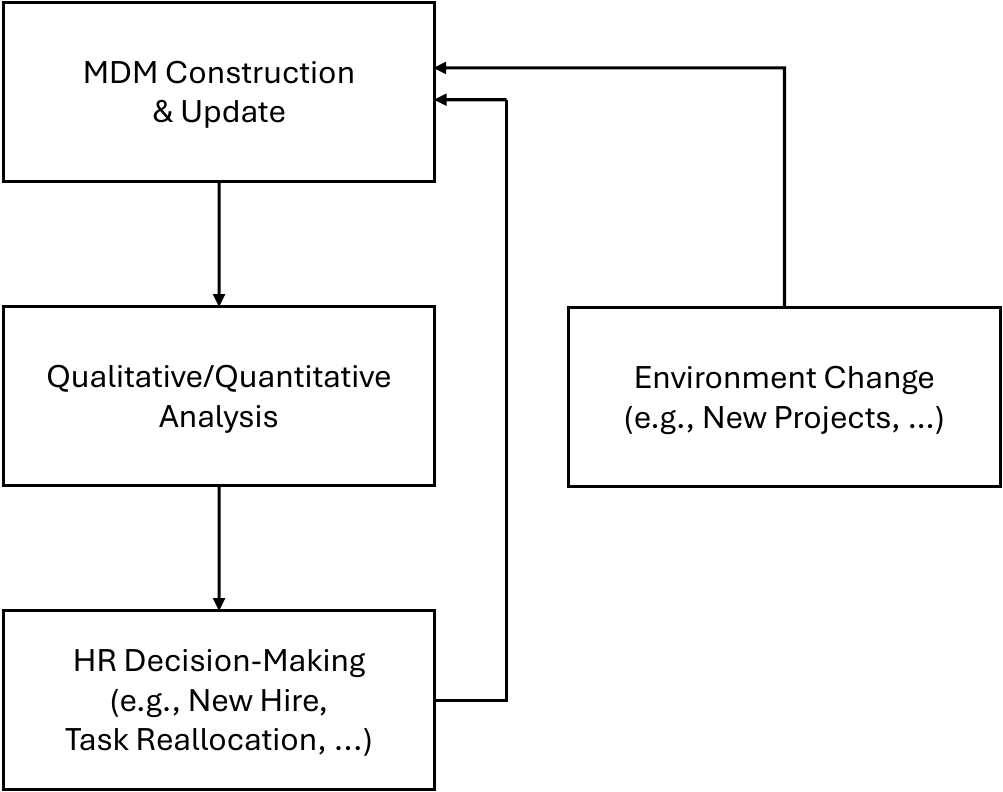}
    \caption{Overall process of the MDM-based HR decision support framework.}
    \label{fig:framework-flow}
\end{figure}

This section presents an MDM-based framework for organizational analysis and HR decision support, with particular relevance to small organizations characterized by complex internal interrelationships. However, the framework is not limited to startups. It is formulated in generalizable terms and can therefore be applied to organizations of varying sizes and across diverse industry contexts, with specific instantiation details determined by practitioners according to organizational characteristics. Figure~\ref{fig:framework-flow} illustrates the overall framework process, which follows a core workflow from initial MDM construction through qualitative/quantitative analysis to HR decision-making. The framework also incorporates an MDM update feedback loop that enables iterative refinement: analysis results inform revisions to the MDM structure, while HR decisions trigger reassessment of the organizational model. In addition, changes in the external environment (e.g., new projects or emerging technologies) feed into the MDM update process, ensuring that the framework remains responsive to evolving conditions.

The remainder of this section describes each module in detail.

\subsection{Initial MDM Construction}
\label{subsec:mdm-construction}

\subsubsection{Domain and Interface Definitions}
\label{subsubsec:domain-interface}

\begin{figure*}[!t]
    \centering

    \includegraphics[page=4, width=0.8\textwidth]{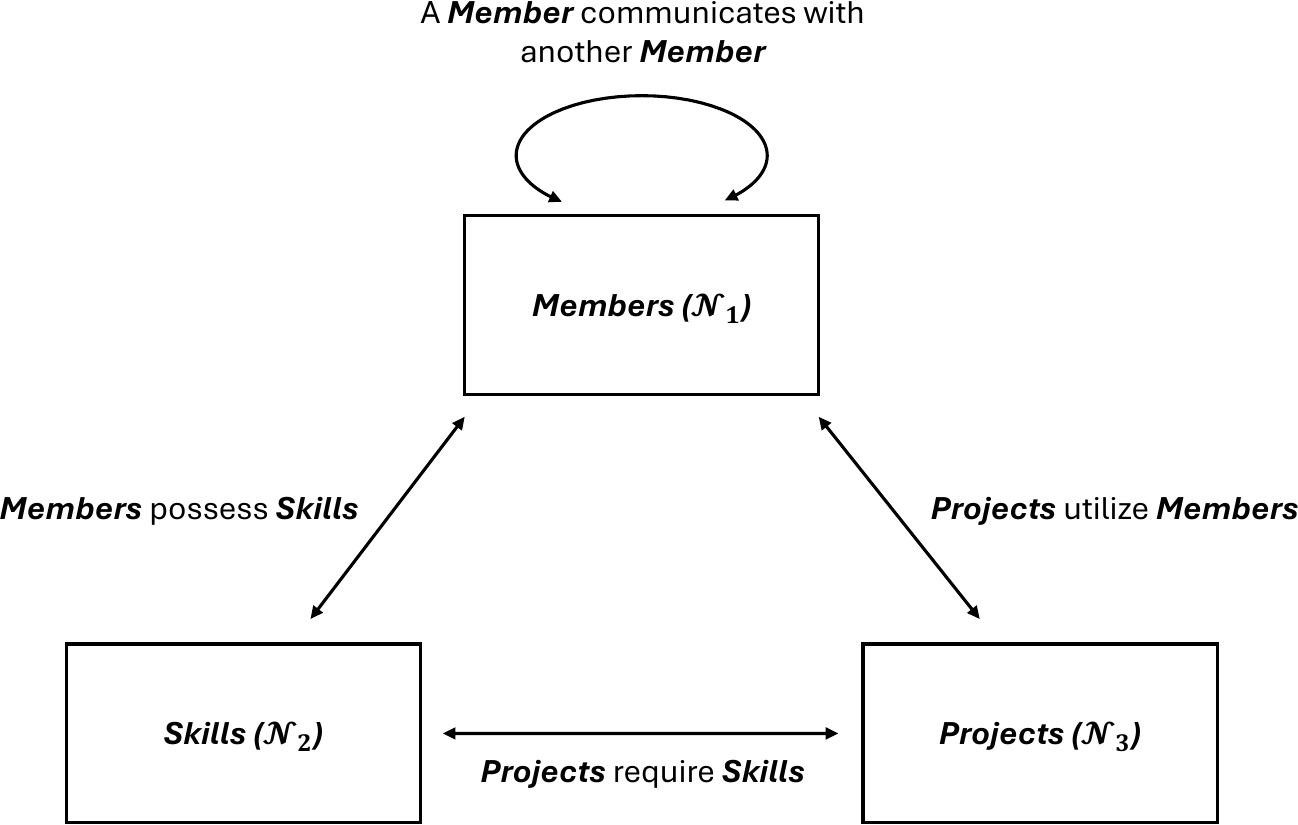}
    \caption{Domain interactions among Members, Skills, and Projects.}
    \label{fig:domain-interactions-diagram}

    \vspace{0.8em}

    \captionsetup{type=table}
    \caption{Interfaces, the domain relations they represent, and their scoring schemes.}
    \label{tab:domain-interactions-table}

    \vspace{0.3em}

    \rowcolors{2}{white}{gray!10}
    \footnotesize
    \setlength{\tabcolsep}{6pt}
    \begin{tabular}{llll}
    \rowcolor{gray!25}
    \textbf{Interface} & \textbf{Relation} & \textbf{Domains} & \textbf{Score} \\
    Work Connectivity & Members communicate with Members & $\mathcal{N}_1$--$\mathcal{N}_1$ & Frequency \\
    Capability Distribution & Members possess Skills & $\mathcal{N}_1$--$\mathcal{N}_2$ & Level of Expertise \\
    Participation Patterns & Projects utilize Members & $\mathcal{N}_1$--$\mathcal{N}_3$ & Level of Contribution \\
    Skill Utilization & Projects require Skills & $\mathcal{N}_2$--$\mathcal{N}_3$ & Binary \\
    \end{tabular}

\end{figure*}

The framework employs three domains that capture the essential dimensions of startup organizational structure: \textit{Members}, \textit{Skills}, and \textit{Projects}. The Members domain ($\mathcal{N}_1$) comprises the set of organizational members, including their roles and responsibilities, and should encompass all personnel whose workload and contributions are relevant to HR decision-making. The Skills domain ($\mathcal{N}_2$) captures the technical and professional competencies required for organizational activities and should be defined at a level of granularity that enables meaningful differentiation among team members while remaining tractable for data collection. The Projects domain ($\mathcal{N}_3$) represents the set of ongoing or planned projects, initiatives, or work streams, defined such that they correspond to distinct work efforts with identifiable participants and required skills. Figure~\ref{fig:domain-interactions-diagram} illustrates the interactions among these three domains through defined interfaces.




Specifically, the framework focuses on four key interfaces---Members--Members, Members--Skills, Members--Projects, and Skills--Projects---that are most directly relevant to HR decision-making. The Members--Members interface captures communication frequency between members (e.g., daily or weekly). The Members--Skills interface represents each member's proficiency in each skill (e.g., based on prior experience or assessed expertise level). The Members--Projects interface indicates each member's level of contribution to each project (e.g., participant or project lead). The Skills--Projects interface specifies the skills required for each project, together with the required proficiency level for each skill (e.g., practitioner or expert). Figure~\ref{fig:domain-interactions-diagram} and Table~\ref{tab:domain-interactions-table} summarize these four interfaces, the domain relations they represent, and their associated scoring schemes.


\subsubsection{Construction Process}
\label{subsubsec:construction-process}
The MDM construction process proceeds in three stages: domain element identification, interface measurement, and matrix assembly and validation. First, the relevant elements in each domain are identified and enumerated. The member domain includes all team members whose work relationships are to be analyzed; in the Skills domain, the core competencies required for organizational activities are defined at an appropriate granularity; and in the Projects domain, the set of ongoing and planned work efforts is specified. Next, systematic data collection is conducted for each interface based on the criteria introduced earlier. The Members--Members interface can be measured through communication records, such as messaging-platform data, email exchanges, and meeting logs. The Members--Skills interface is typically generated through self-assessment and peer evaluation, while the Members--Projects interface is determined from project management tools and work logs. The Skills--Projects interface identifies the skills required for each project. Finally, all collected data are quantified according to the defined scales and assembled into the MDM, after which the matrix is validated for consistency and accuracy. Cross-validation by multiple reviewers is recommended to improve the reliability of the resulting representation.

\begin{figure}[h!]
    \centering
        \includegraphics[page=2, width=0.9\columnwidth]{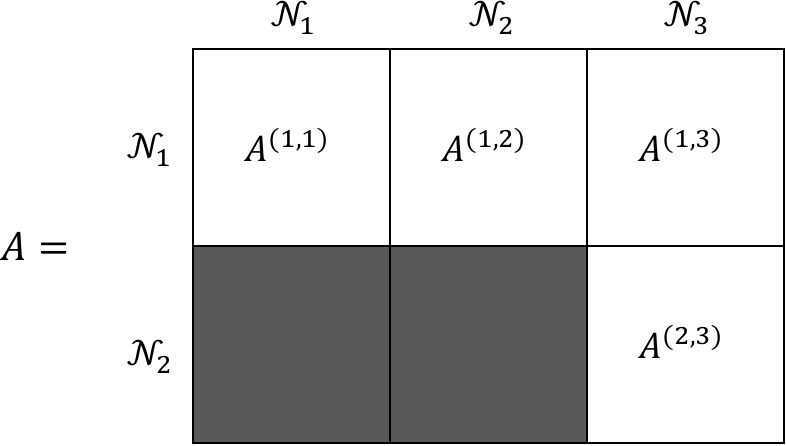}
    \caption{Proposed MDM structure with one DSM ($A^{(1,1)}$) and three DMMs ($A^{(1,2)}$, $A^{(1,3)}$, and $A^{(2,3)}$).}
    \label{fig:matrix-focus}
\end{figure}

The resulting MDM is denoted by $A = [a_{ij}]$, where $a_{ij}$ represents the interface value between elements $i$ and $j$ where $i$ and $j$ belong to different domains. Figure~\ref{fig:matrix-focus} illustrates the structure of the MDM used in this study, which comprises one DSM block and three DMM blocks. Each block can be written in the general form
\begin{equation}
A^{(d_1,d_2)} = [a_{ij}]_{i \in \mathcal{N}_{d_1},\, j \in \mathcal{N}_{d_2}},
\label{eq:submatrix}
\end{equation}
where $A^{(d_1,d_2)}$ denotes the DSM or DMM capturing the relationships between domains $\mathcal{N}_{d_1}$ and $\mathcal{N}_{d_2}$. In this framework, the relevant blocks are Members $\times$ Members ($A^{(1,1)}$), Members $\times$ Skills ($A^{(1,2)}$), Members $\times$ Projects ($A^{(1,3)}$), and Skills $\times$ Projects ($A^{(2,3)}$).

\subsection{MDM Construction}
\label{subsec:qualitative-analysis}
Qualitative analysis of the constructed MDM enables identification of organizational patterns and anomalies through visual inspection of each interface region. This section provides interpretation guidelines for the four key interfaces.

\subsubsection{Members--Members: Communication Patterns}
\label{subsubsec:comm-patterns}

The Members--Members interface reveals the communication structure of the organization. Several recurring patterns are particularly noteworthy. Hub--spoke structures arise when certain members maintain high communication frequency with many others and thereby function as information hubs. While such structures can support coordination efficiency, excessive dependence on specific individuals may create bottlenecks and single points of failure. Communication silos---groups of members characterized by dense internal communication but limited external interaction---may indicate team isolation or coordination gaps. Bridge members, who connect otherwise disconnected groups, play a critical role in facilitating cross-functional collaboration.

\subsubsection{Members--Skills: Capability Distribution}
\label{subsubsec:skill-distribution}

The Members--Skills interface maps organizational capabilities to personnel. Skill concentration---where only one member possesses a particular skill at the expert level---creates a single point of failure and thus an organizational risk. In contrast, skill redundancy---where multiple members possess expert-level proficiency in the same skill---improves resilience, although it may also suggest overinvestment. Capability gaps, namely skills lacking expert-level coverage, may indicate the need for targeted training or hiring.

\subsubsection{Members--Projects: Participation Patterns}
\label{subsubsec:participation-patterns}

The Members--Projects interface reveals how workload is distributed across the organization. Load concentration arises when certain members serve as primary contributors on multiple projects at the same time, increasing their risk of overload and burnout. Underutilization, where members have limited project involvement, may indicate role misalignment or inefficient capacity allocation. Project dependency---where projects rely heavily on specific individuals---introduces schedule and quality risks.

\subsubsection{Skills--Projects: Skill Utilization}
\label{subsubsec:skill-utilization}

The Skills--Projects interface identifies the competencies required across the project portfolio. High-demand skills---required by multiple projects---may justify deeper investment in capability development. In contrast, underutilized skills not associated with current projects may suggest strategic misalignment or potential future opportunity areas. Projects requiring skills that are not adequately covered by existing personnel present execution risks due to insufficient skill coverage.

\subsection{Quantitative Analysis}
\label{subsec:quantitative-analysis}

\subsubsection{Base Metrics}
\label{subsubsec:metrics}
The framework defines four base metrics that capture distinct dimensions of organizational structure: \textit{Communication Score}, \textit{Project Role Score}, \textit{Skill Leverage}, and \textit{Skill Gap}. Each metric is derived from the MDM interfaces and is computed for each member in $\mathcal{N}_1$.

The first metric, the Communication Score (Comm), measures the extent to which a team member is engaged in organizational communication. A high Comm value indicates that the member functions as an information hub, while also reflecting the substantial time and effort devoted to communication activities. For member $i$, Comm is defined as
\begin{equation}
\text{Comm}_i = \sum_{j\in\mathcal{N}_1-{\left\{i\right\}}} a_{ij}^{(1,1)},\quad \forall i\in \mathcal{N}_1,
\label{eq:comm}
\end{equation}
where $a_{ij}^{(1,1)}$ denotes the communication-frequency weight between members $i$ and $j$ in the Members--Members DSM $A^{(1,1)}$.

The second metric, the Project Role Score (Proj), quantifies the intensity of a team member's project participation. By assigning different weights to different roles, it reflects the level of responsibility borne by each member rather than merely counting the number of projects in which the member participates. For member $i$, Proj is defined as
\begin{equation}
\text{Proj}_i = \sum_{k\in\mathcal{N}_3} a_{ik}^a{(1,3)},\quad \forall i\in \mathcal{N}_1,
\label{eq:proj}
\end{equation}
where $a_{ik}^{(1,3)}$ denotes the participation weight of member $i$ in project $k$ in the Members--Projects DMM $A^{(1,3)}$.

The third metric, Skill Leverage (SL), measures the extent to which a team member's skills are utilized in areas of organizational demand. Rather than considering skill possession alone, this metric reflects whether a member's skills are associated with competencies that are required across the project portfolio. For member $i$, SL is defined as
\begin{equation}
\text{SL}_i = \sum_{j\in\mathcal{N}_2} \left( a_{ij}^{(1,2)} \times u_j \right),\quad \forall i\in \mathcal{N}_1,
\label{eq:sl}
\end{equation}
where $a_{ij}^{(1,2)}$ denotes the proficiency level of member $i$ in skill $j$ in the Members--Skills DMM $A^{(1,2)}$, and $u_j$ represents the organizational demand for skill $j$, defined as the number of projects that require that skill:
\begin{equation}
u_j=\sum_{k\in\mathcal{N}_3} \mathbf{1}[a_{jk}^{(2,3)} > 0],\quad \forall j\in\mathcal{N}_2.
\label{eq:uj}
\end{equation}

The fourth metric, the Skill Gap (Gap), measures deficiencies in the skills required for the projects to which a member is assigned. A high Gap value implies hidden costs such as learning overhead, mentoring requirements, and reduced work efficiency. Let
\begin{equation}
P_i = \{k \in \mathcal{N}_3 : a_{ik}^{(1,3)} > 0\}
\label{eq:pi}
\end{equation}
denote the set of projects in which member $i$ participates, and let
\begin{equation}
R_i = \{j \in \mathcal{N}_2 : \exists\, k \in P_i,\; a_{jk}^{(2,3)} > 0\}
\label{eq:ri}
\end{equation}
denote the set of skills required by at least one project in $P_i$. Then, for member $i$, Gap is defined as
\begin{equation}
\text{Gap}_i = \frac{1}{|R_i|}\sum_{j\in R_i} \max(0,\, \theta - a_{ij}^{(1,2)}),\quad \forall i\in \mathcal{N}_1,
\label{eq:gap}
\end{equation}
where $\theta$ is a threshold representing the proficiency level considered adequate for effective project execution. The inner $\max(0,\cdot)$ operator ensures that only skill deficiencies contribute to the Gap, making the metric non-compensatory: surplus proficiency in one skill cannot offset a deficiency in another, since the learning overhead and mentoring cost of an under-covered skill are incurred independently of a member's expertise in other areas.

To facilitate the evaluation of each member's metric values within a given organizational state and support relative comparison across different metrics, the basic metrics are reported using both absolute and normalized values. The normalized values are obtained using min--max scaling, which maps each metric to the interval $[0,1]$, as follows:
\begin{equation}
\hat{s}_i = \frac{s_i - s_{\min}}{s_{\max} - s_{\min}},
\label{eq:minmax}
\end{equation}
where $s_{\min}=\min_{i\in \mathcal{N}_1}s_i$ and $s_{\max}=\max_{i\in \mathcal{N}_1}s_i$ denote the minimum and maximum values of the corresponding metric across all members, respectively.

\subsubsection{Composite Indices and Quadrant Analysis}
\label{subsubsec:quadrant-analysis}

The base metrics are further combined into two composite indices that synthesize the workload and value dimensions of each team member: the \textit{Workload Index} and the \textit{Value Index}.

The Workload Index (WI) measures the overall work burden borne by each team member. For member $i$, WI is defined as
\begin{equation}
\text{WI}_i = w_1 \cdot \widehat{\text{Proj}}_i + w_2 \cdot \widehat{\text{Comm}}_i + w_3 \cdot \widehat{\text{Gap}}_i,
\label{eq:workload}
\end{equation}
where $w_1$, $w_2$, and $w_3$ are non-negative weight parameters satisfying $w_1 + w_2 + w_3 = 1$. These weights should be specified according to the organizational context. As a weighted aggregation of normalized Proj, Comm, and Gap values, the WI captures the extent to which each member is burdened by project participation, communication demands, and learning requirements associated with assigned projects.

The Value Index (VI) quantifies the organizational value contributed by each team member. For member $i$, VI is defined as
\begin{equation}
\text{VI}_i = v_1 \cdot \widehat{\text{SL}}_i + v_2 \cdot \widehat{\text{Proj}}_i + v_3 \cdot \widehat{\text{Comm}}_i,
\label{eq:value}
\end{equation}
where $v_1$, $v_2$, and $v_3$ are non-negative weight parameters satisfying $v_1 + v_2 + v_3 = 1$. Similar to the WI, these weights should be specified according to organizational priorities. As a weighted aggregation of normalized SL, Proj, and Comm values, the VI captures each member's organizational contribution in terms of skill level, project participation, and communication involvement.

The weight values may be established through stakeholder input, such as interviews with team leads or executive judgment, or through more systematic approaches such as the analytic hierarchy process (AHP)~\cite{saaty1977scaling,saaty1980ahp}. The framework does not prescribe fixed weight values, as appropriate choices depend on organizational characteristics and strategic priorities.

\begin{figure}
    \centering
        \includegraphics[page=5, width=0.99\columnwidth]{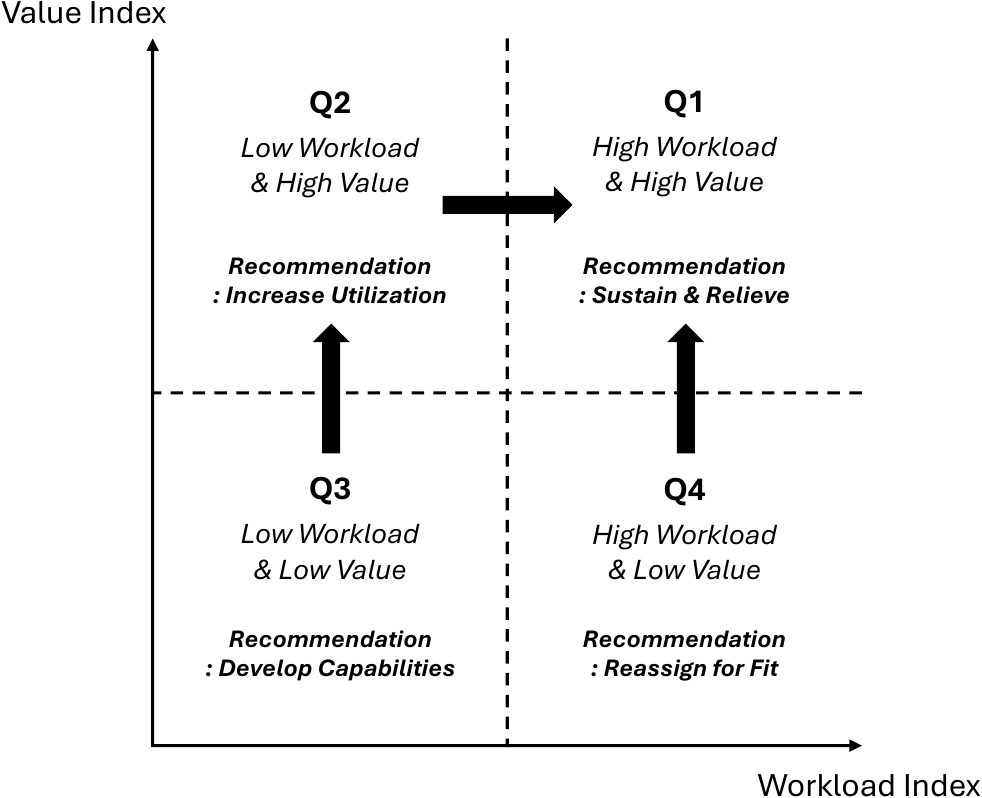}
    \caption{Quadrant analysis framework with HR recommendations. Arrows indicate the desired movement in the workload-value space that each recommendation aims to achieve.}
    \label{fig:quadrant-framework}
\end{figure}

The WI and VI are used to position team members in a two-dimensional space (Figure~\ref{fig:quadrant-framework}), with quadrant boundaries defined by the mean values of the two indices. Using the mean as the boundary yields a relative classification of each member's workload and value contribution within the current team composition. This relative framing enables two complementary directions for HR decision-making at any given point in time.

The first direction concerns quadrant membership: members outside Q1 represent opportunities for development, and HR interventions---such as workload assignment, skill development, or project reallocation---should aim to progressively move them toward Q1, where both engagement and value contribution are high. The second direction concerns extreme positioning within Q1: a member with WI approaching 1.0 carries a disproportionate share of the team's workload relative to peers, making their situation unsustainable even within the ideal quadrant. In such cases, the appropriate intervention is workload redistribution to reduce that member's WI. Here, the mean boundary serves as a guardrail: the goal is to bring extreme WI values back toward the center, but not below the mean, as crossing that threshold would displace the member from Q1 entirely and signal underutilization rather than relief.

Together, these two directions define the decision logic of the framework: lift underperforming members toward Q1, and relieve overburdened Q1 members while keeping them within it. This representation partitions members into four categories, each associated with distinct HR implications and recommended actions.

The resulting quadrant interpretation provides actionable guidance for HR decision-making. Q1 (High Workload / High Value) represents the state in which high-value members handle a large volume of work. In startups, where maintaining a high density of capable talent is essential, this is arguably the ideal operational state. However, sustained high workload carries burnout risk; therefore, the recommended actions include \textit{appropriate compensation}, provision of adequate rest, or \textit{hiring} to share the workload burden. Q2 (Low Workload / High Value) indicates that high-value members are not being fully utilized. The recommended action is to \textit{increase utilization} by assigning additional responsibilities that leverage their capabilities. Q3 (Low Workload / Low Value) captures members whose current value contribution is low and who consequently have limited task assignments. The recommended action is \textit{capability development} through training and skill-building, laying the groundwork for these members to deliver higher value and take on greater responsibilities in the future. Finally, Q4 (High Workload / Low Value) suggests that members are carrying heavy workloads on tasks that may not align well with their skill sets. The recommended action is \textit{reassignment} to projects with better skill fit, which can increase their value contribution without necessarily reducing their workload.

By systematically applying these quadrant-specific interventions---sustaining Q1 members, leveraging Q2 members, developing Q3 members, and realigning Q4 members---the organization progressively raises both its mean WI and mean VI over time. As these mean values shift upward, the quadrant boundaries shift accordingly, setting a higher standard for what constitutes above-average performance. The quadrant boundaries adapt as the organization improves, ensuring that the analysis remains a driver of continuous advancement rather than a static snapshot. In this way, the quadrant framework serves not only as a diagnostic tool for prioritizing immediate HR interventions but also as a mechanism for sustained organizational growth.

\section{Case Study: Hiring Decision at Planby Technologies, Inc.}
\label{sec:case-study}

To demonstrate the framework's practical applicability and provide insights for practitioners in similar startup environments, we applied the proposed framework to an early-stage technology startup, Planby Technologies, Inc. (\url{https://planby.us}).

\begin{figure}[h!]
    \centering
    \includegraphics[width=0.98\columnwidth]{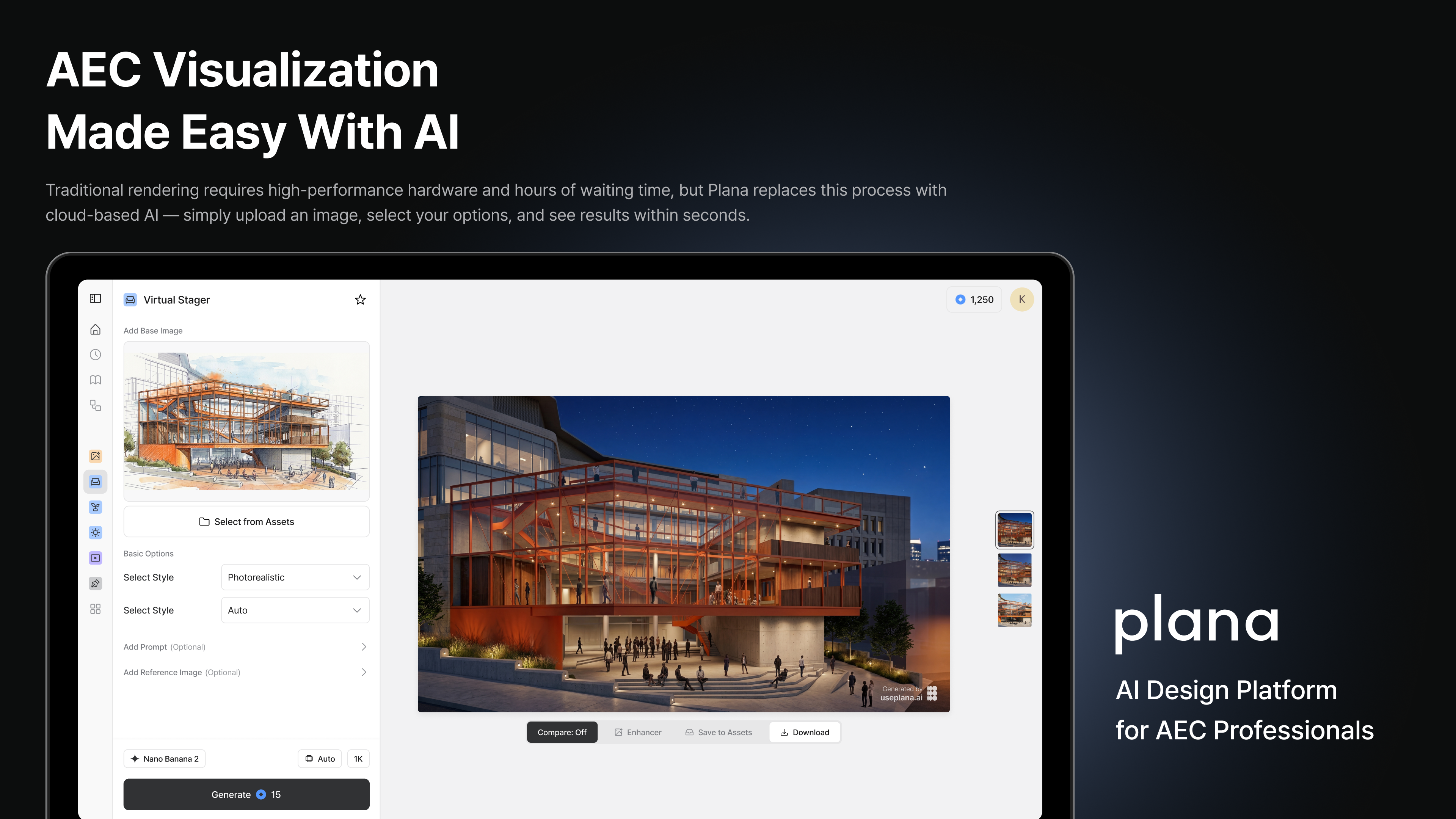}
    \caption{Representative interface of Plana, Planby Technologies' AI design platform for AEC professionals.}
    \label{fig:planby-product}
\end{figure}

Planby Technologies, Inc. is an early-stage startup that builds AI-powered visualization tools for the architecture, engineering, and construction (AEC) industry. Its flagship product, Plana (\url{https://useplana.ai/}), allows AEC professionals to transform conceptual sketches and base images into photorealistic renderings within seconds through a cloud-based AI workflow, replacing traditional rendering pipelines that typically require high-performance hardware and hours of waiting time. A representative view of the product interface is shown in Figure~\ref{fig:planby-product}. The organization comprises 13 team members ($|\mathcal{N}_1| = 13$) with diverse backgrounds spanning software engineering, AI/ML research, product design, and business development. As a rapidly growing startup, Planby Technologies, Inc. faces the typical challenges of resource allocation, skill coverage, and workload distribution that characterize early-stage technology ventures.

The case study is structured as a two-phase analysis. In the first phase, we construct the MDM from a three-month observation period and perform qualitative and quantitative assessments of the current organizational state. In the second phase, we introduce a strategic hiring decision informed by the first-phase findings and evaluate its projected impact on workload distribution and organizational value.

\subsection{Phase 1: Current State Analysis}
\label{subsec:current-state}

\subsubsection{Initial MDM Construction}
\label{subsubsec:initial-mdm}

Following the framework guidelines (Section~\ref{sec:framework}), we instantiated the three domains as follows. The Members domain ($|\mathcal{N}_1| = 13$) includes all current team members: developers, designers, project managers, and account managers. The Skills domain ($|\mathcal{N}_2| = 10$) comprises ten core competencies identified as critical for operations: Product Management, UI/UX Design, Frontend Development, Backend Development, LLM Engineering, Vision Engineering, Business Development, Sales, Marketing, and Operations. The Projects domain ($|\mathcal{N}_3| = 6$) encompasses six ongoing projects, including new product development, existing service improvement, customer support, and deliverable maintenance initiatives.

\begin{figure*}
    \centering
        \includegraphics[page=6, width=0.99\textwidth]{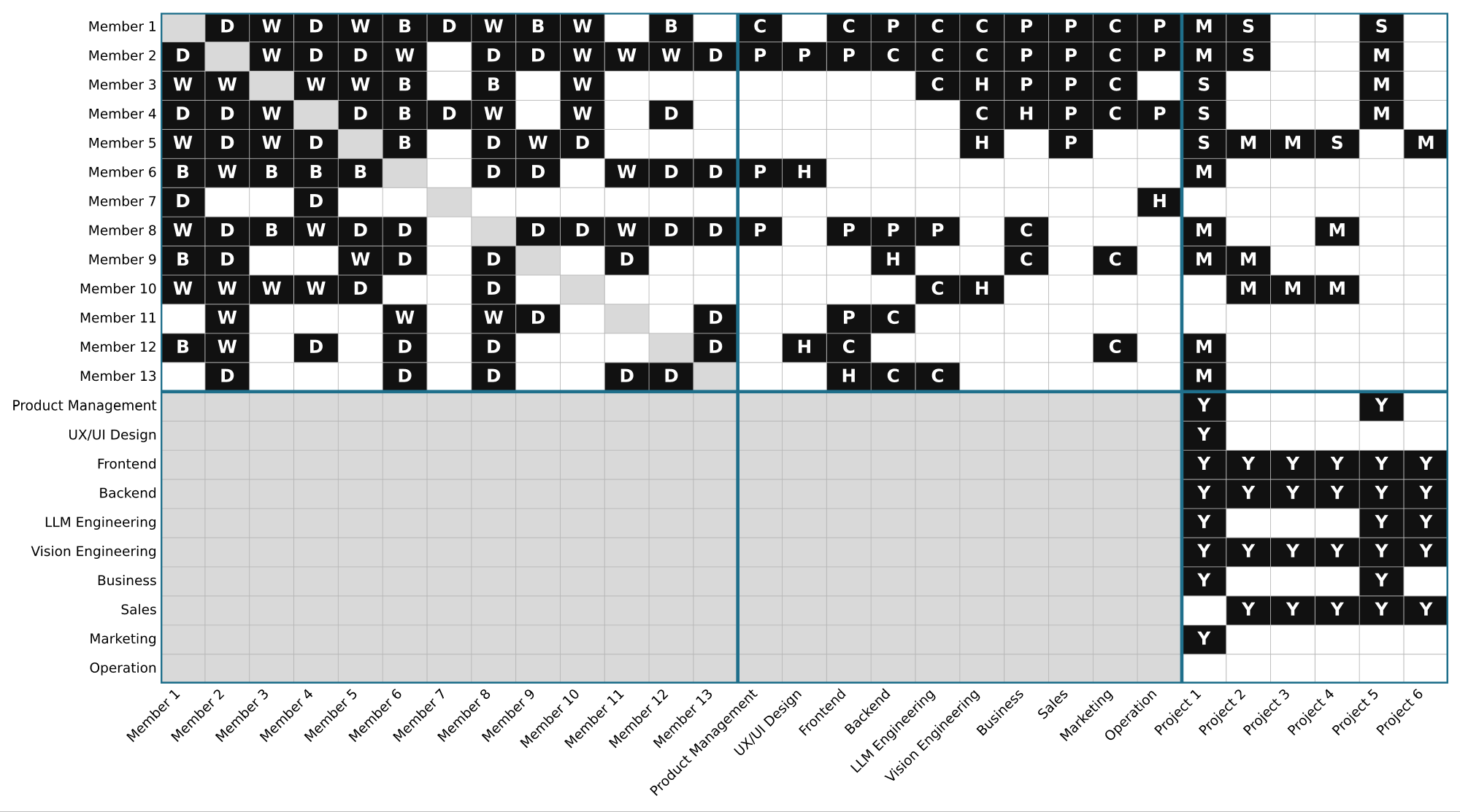}
    \caption{MDM of Planby Technologies, Inc. for Phase 1.}
    \label{fig:planby-mdm}
\end{figure*}

Data was collected over a three-month observation period. For the Members--Members interface, communication frequency was measured through analysis of Slack messages, email exchanges, and meeting attendance records, with frequencies classified as Daily (D, weight 3.0), Weekly (W, weight 2.0), or Bi-weekly (B, weight 1.0). For Members--Skills, proficiency levels were assessed through self-evaluation and peer review, categorized as Expert (H, weight 5.0), Practical (P, weight 3.0), or Conceptual (C, weight 1.0). Members--Projects participation roles were determined from project management tools and verified through project leader interviews, with Main (M, weight 4.0) and Sub (S, weight 1.0) designations. Finally, Skills--Projects skill utilization was identified from project technology stacks. Table~\ref{table:interfaces} summarizes the data sources, classification levels, and assigned weights for each MDM interface. The constructed MDM is shown in Figure~\ref{fig:planby-mdm}.

\begin{table}[h!]
\centering
\caption{Summary of MDM interface types, classification levels, and assigned weights.}
\label{table:interfaces}
\rowcolors{2}{white}{gray!10}
\scriptsize
\setlength{\tabcolsep}{3pt}
\begin{tabular}{llll}
\rowcolor{gray!25}
\textbf{Interface} & \textbf{Data Source} & \textbf{Classification} & \textbf{Weight} \\

\multirow[t]{3}{*}{Members--Members} 
& \multirow[t]{3}{*}{Slack, email, meetings} 
& Daily (D) & 3.0 \\
& & Weekly (W) & 2.0 \\
& & Bi-weekly (B) & 1.0 \\

\multirow[t]{3}{*}{Members--Skills} 
& \multirow[t]{3}{*}{Self-evaluation, peer review}
& Expert (H) & 5.0 \\
& & Practical (P) & 3.0 \\
& & Conceptual (C) & 1.0 \\

\multirow[t]{2}{*}{Members--Projects}
& \multirow[t]{2}{*}{PM tools, interviews}
& Main (M) & 4.0 \\
& & Sub (S) & 1.0 \\

\multirow[t]{2}{*}{Skills--Projects}
& \multirow[t]{2}{*}{Project tech stacks}
& Required & 1.0 \\
& & Not required & 0.0 \\
\end{tabular}
\end{table}

\subsubsection{Qualitative Analysis}
\label{subsubsec:qualitative-results}

Analysis of the Members--Members interface revealed that organizational communication follows a Hub-Spoke structure centered around select individuals. Notably, Members 2 and 8 exhibited exceptionally high communication frequencies, confirming their roles as primary conduits for information flow. Member 8, in particular, communicates daily with 7 out of the 13 team members, suggesting a pivotal, cross-functional role.

In contrast, certain team members showed relatively low connectivity, often due to their isolated focus on specialized tasks. This communication imbalance can lead to information silos and place an undue burden on core members. Visual inspection also highlighted a communication gap between the business and product teams, while intra-team communication among developers (Members 8, 9, 12 and 13) and designers (Members 6 and 12) remained robust.

In terms of skill distribution, while overall expertise is high, the matrix identified insufficient redundancy in mission-critical skills. For instance, only one member possesses expert-level (H) competency in each of Backend Development and Operations (Members 7 and 9), creating a significant single-point-of-failure risk.

Conversely, areas like Frontend Development and Sales are well-supported, with multiple members possessing practical (P) level proficiency or higher. However, while Backend Development shows a skill distribution similar to Frontend Development, a key distinction lies in the availability of those skills for project work. Members 1 and 8, who hold Practical-level or higher competency in Backend Development, also possess a broad range of other skills and are consequently allocated to projects requiring those diverse capabilities. As a result, Member 9 effectively shoulders the Backend Development responsibilities across projects largely alone, representing a concentrated dependency on a single individual. While Vision Engineering has three expert-level practitioners, these are concentrated among already high-workload members (Members 5 and 10), creating a redundancy risk if key personnel depart. AI/ML skills showed high interdependency with core development skills, emphasizing the necessity for an integrated technical capability development plan.

Regarding project participation, the analysis revealed a significant concentration of responsibility. Members 5 and 10 are currently serving as primary personnel on 3 out of 6 active projects, resulting in the highest measured project loads. Members 8 and 9 also maintain high work intensities, managing primary roles across 2 projects each. This concentration stems from (1) the high demand for advanced technical specialized skills, (2) the technical complexity of specific projects, and (3) a lack of experienced mid-level personnel to share the load.

Notably, Members 7 and 11 have no project assignments in the current MDM. Member 7's sole expert-level skill is Operations, which supports internal organizational functions such as financial management and operational processes rather than project deliverables---explaining its absence from the Skills--Projects matrix. Member 11, while possessing Frontend (P) and Backend (C) proficiency, similarly lacks project participation. These two members represent 15\% of the team operating outside the project structure. The framework surfaces this through low Workload and Value indices, though practitioners should note that internally focused roles like Operations contribute organizational value not fully captured by project-based metrics.

\subsubsection{Quantitative Analysis}
\label{subsubsec:quantitative-results}
To complement the qualitative observations, we computed the base metrics and composite indices defined in Section~\ref{subsec:quantitative-analysis} for the Planby Technologies team. For the base metrics, the Comm, SL, and Proj were computed directly from the interface weights described above. For the Gap, the proficiency threshold was set to $\theta = 3.0$, corresponding to the Practical level, meaning that members holding only conceptual level proficiency in a project-required skill are considered to have a gap.

Based on interviews with team members and leadership, the following composite index weights were adopted: for the WI, $w_1 = 0.6$ (Proj), $w_2 = 0.3$ (Comm), $w_3 = 0.1$ (Gap); for the VI, $v_1 = 0.5$ (SL), $v_2 = 0.3$ (Proj), $v_3 = 0.2$ (Comm). These weights reflect the organization's project-intensive nature and emphasis on technical expertise as a primary value driver.
Using min-max normalization (Equation~\ref{eq:minmax}), the four base metrics and two composite indices were computed for all 13 team members, as presented in Table~\ref{table:scores}. In the table, each normalized value is accompanied by its raw score in parentheses, and the composite indices (WI and VI) include their rank among all members in parentheses.

\begin{table}[h!]
\centering
\caption{Normalized base metrics and composite indices for all 13 team members in Phase 1. Original values are shown in parentheses for the basic metrics, whereas rankings are shown in parentheses for the composite indices.}
\label{table:scores}
\rowcolors{2}{white}{gray!10}
\scriptsize
\setlength{\tabcolsep}{3pt}
\begin{tabular}{l|cccc|cc}
\rowcolor{gray!25}
\textbf{Mem.} & \textbf{Comm} & \textbf{SL} & \textbf{Proj} & \textbf{Gap} & \textbf{WI} & \textbf{VI} \\
1  & 0.64 (20) & 0.89 (57) & 0.43 (6)  & 0.59 (1.44) & 0.51 (7)  & 0.70 (4) \\
2  & 1.00 (28) & 1.00 (64) & 0.64 (9)  & 0.36 (0.89) & 0.72 (2)  & 0.89 (1) \\
3  & 0.27 (12) & 0.86 (55) & 0.36 (5)  & 0.73 (1.78) & 0.37 (11) & 0.59 (5) \\
4  & 0.73 (22) & 0.50 (32) & 0.36 (5)  & 0.86 (2.11) & 0.52 (6)  & 0.50 (8) \\
5  & 0.59 (19) & 0.70 (45) & 1.00 (14) & 0.95 (2.33) & 0.87 (1)  & 0.77 (3) \\
6  & 0.64 (20) & 0.17 (11) & 0.29 (4)  & 0.92 (2.25) & 0.45 (8)  & 0.30 (10) \\
7  & 0.00 (6)  & 0.00 (0)  & 0.00 (0)  & 0.00 (0.00) & 0.00 (13) & 0.00 (13) \\
8  & 1.00 (28) & 0.83 (53) & 0.57 (8)  & 0.64 (1.56) & 0.71 (4)  & 0.79 (2) \\
9  & 0.41 (15) & 0.52 (33) & 0.57 (8)  & 1.00 (2.44) & 0.57 (5)  & 0.51 (7) \\
10 & 0.36 (14) & 0.52 (33) & 0.86 (12) & 0.92 (2.25) & 0.72 (3)  & 0.59 (6) \\
11 & 0.27 (12) & 0.38 (24) & 0.00 (0)  & 0.00 (0.00) & 0.08 (12) & 0.24 (12) \\
12 & 0.41 (15) & 0.19 (12) & 0.29 (4)  & 0.97 (2.38) & 0.39 (9)  & 0.26 (11) \\
13 & 0.41 (15) & 0.61 (39) & 0.29 (4)  & 0.97 (2.38) & 0.39 (9)  & 0.47 (9) \\
\end{tabular}
\end{table}

\begin{figure}[h!]
\centering
\includegraphics[page=10, width=0.99\columnwidth]{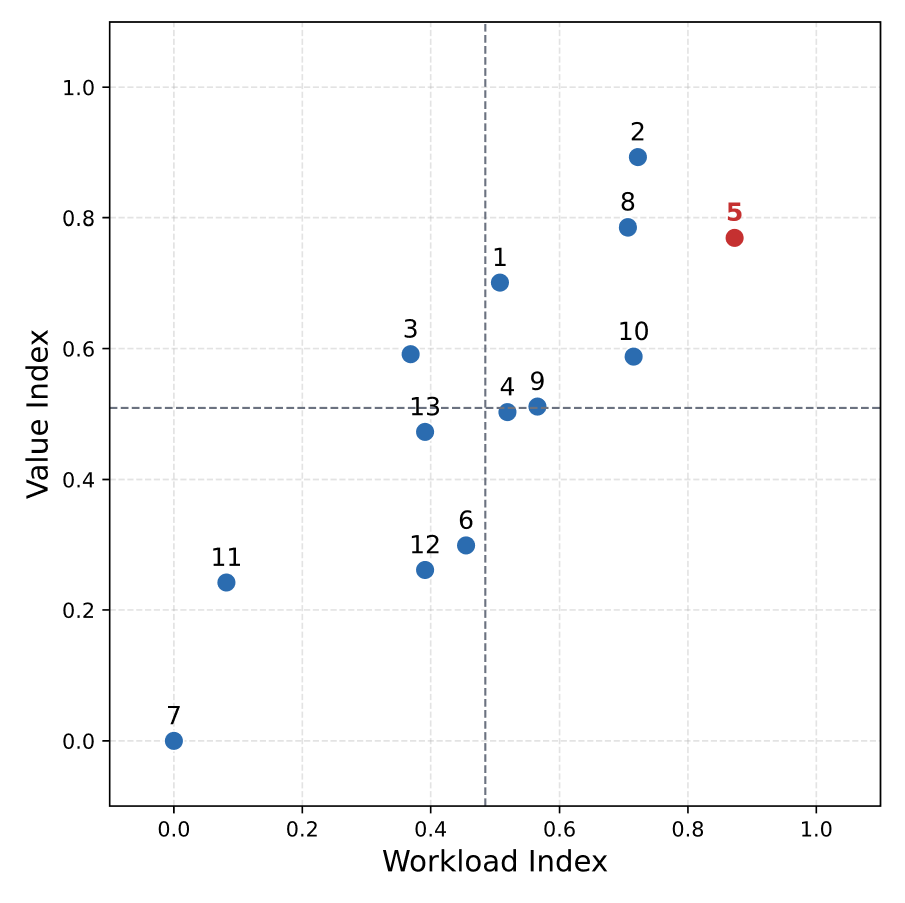}
\caption{Quadrant analysis of Phase 1.}
\label{fig:quadrant-pre}
\end{figure}

Section~\ref{subsubsec:quadrant-analysis} prescribes the HR action appropriate to each quadrant: hiring or workload redistribution for Q1, expanded responsibilities for Q2, capability development for Q3, and project reassignment for Q4. Figure~\ref{fig:quadrant-pre} positions the 13 Planby members in the WI--VI plane, and the difficulty of the resulting decision within each quadrant depends on how many members it contains.

For sparsely populated quadrants---Q2 with Member 3 alone and Q4 with Member 4 alone---the prescribed action applies in a straightforward manner: each single member is an unambiguous target for the corresponding intervention, and qualitative analysis together with routine organizational awareness already suffices to identify them. The quantitative indices add little beyond confirmation in such cases.

The decision problem becomes substantive when many members fall into the same quadrant, as in the Q1 cluster of six (Members 1, 2, 5, 8, 9, and 10) or the Q3 cluster of five (Members 6, 7, 11, 12, and 13). The prescribed action cannot be applied uniformly to every member within the budget and bandwidth of a startup, and the operative question shifts from \emph{what to do} to \emph{to whom first}. Focusing on Q1, where workload relief is the most time-critical concern, direct observation cannot reliably rank the six members: all appear busy and valuable from a qualitative standpoint, and even experienced leadership may struggle to distinguish their relative urgency. The WI ordering resolves this ambiguity. As Section~\ref{subsubsec:quadrant-analysis} anticipates, an extreme WI position within Q1---approaching 1.0---signals unsustainability even within this otherwise ideal quadrant. Member 5's WI of 0.87 separates from the next-highest Q1 member (Member 2 at 0.72) by a 0.15 margin, identifying Member 5 as the unambiguous priority target for the Q1 intervention. This selectivity---naming \emph{who first} and quantifying \emph{by how much}---is the distinctive contribution of the quantitative layer, and it is what makes the targeted hiring decision in Section~\ref{subsec:hiring-decision} defensible rather than discretionary. More generally, the marginal value of the indices over experienced judgment scales with cluster size: small clusters are navigated by intuition, while large, internally similar clusters can only be triaged through explicit, comparable scores.

Notably, the quadrant plot also exhibits a general upward-right trend, in which members with higher VI also tend to carry higher WI. This pattern partly reflects the metric formulation, as project participation contributes positively to both indices. At the same time, it mirrors a characteristic tendency in startups: high-value members attract a disproportionate share of responsibilities, resulting in a natural correlation between contribution and workload.

\subsubsection{Results Summary and Discussion}
\label{subsubsec:pre-hiring-discussion}
Through the current state analysis, the following key problems were identified. First, communication and work within the organization are concentrated on a small number of key personnel, with Member 5 bearing the highest workload burden ($\mathrm{WI}=0.87$). Second, backup personnel for advanced technical skills such as Backend and Vision Engineering are insufficient, acting as risk factors. Third, project participation imbalance has intensified, with some team members handling 3 main projects simultaneously while others participate in only 1--2 projects.

As a result of quadrant analysis, Members 1, 2, 5, 8, 9 and 10 are positioned in Q1 (High Workload / High Value), indicating that immediate burden reduction is needed. Meanwhile, only Member 3 is positioned in Q2 (Low Workload / High Value), suggesting that opportunities for efficient utilization of high-value personnel are limited.

Synthesizing these analysis results, it is clear that reducing the workload of key personnel and securing backup capabilities for specialized skills are the top priorities for sustainable organizational growth. In particular, hiring personnel with a technical skill set similar to Member 5 was judged to be the most effective solution.

\subsection{Phase 2: Hiring Decision and Post-Hiring Analysis}
\label{subsec:hiring-decision}

Based on the findings from Section~\ref{subsec:current-state}, the organization decided to recruit a new member with a skill profile similar to Member 5's, in order to relieve the most acute workload concentration and provide redundancy for the advanced technical skills currently dependent on a single contributor. The hiring criteria prioritized capabilities in the technical skills critical to Member 5's overloaded projects, together with sufficient breadth to absorb adjacent responsibilities. Member 14 was recruited with proficiency in Vision Engineering at Practical (P) level, LLM Engineering at Conceptual (C) level, and Sales at Practical (P) level.

\subsubsection{Hiring Scenario Description}
\label{subsubsec:hiring-scenario}

The primary objective of hiring Member 14 was to relieve Member 5's workload. After a one-week onboarding period, Member 14 was assigned to take over Project~3 (P3), one of the three projects in which Member 5 was participating as a main role. P3 was selected as the entry project for three reasons: (i) it was an actively progressing engagement that would expose Member 14 to the team's working pace, (ii) Member 5 was not the customer-facing point of contact for P3, allowing Member 14 to learn Planby's project execution style without immediate exposure to client-side communication, and (iii) the technical scope aligned closely with Member 14's skill set. Rather than withdrawing from P3 immediately, Member 5 was retained as a sub-role contributor on the same project to facilitate handover.

This setup motivates analysis at two time points. The first captures the state shortly after hiring within the first one to four weeks, in which Member 14 has joined the team and the formal handover on P3 has been initiated but full integration has not yet occurred. The second captures the projected state three months after hiring, by which point Member 14 is assumed to have completed onboarding and to be operating at expected capacity. While the hiring was justified by Member 5's overload, expecting an immediate and full adaptation from a new hire is unrealistic, and as a fast-moving startup Planby could not afford an onboarding window longer than approximately one month. Comparing these two snapshots therefore separates the transient cost of onboarding from the longer-term benefit of redistribution.

\begin{figure*}[h!]
\centering
\begin{subfigure}{\linewidth}
\includegraphics[page=7, width=0.99\columnwidth]{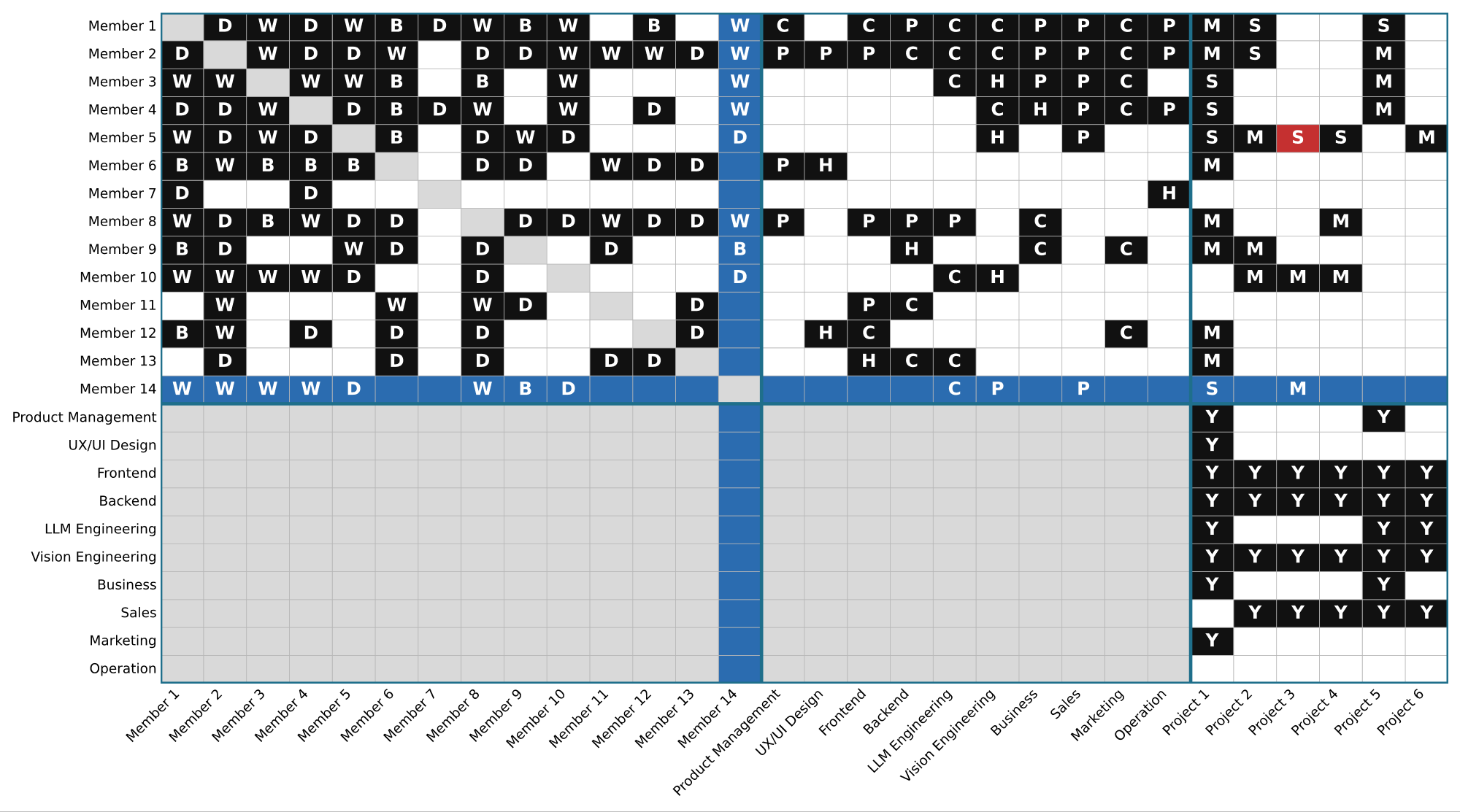}
\caption{Updated MDM immediately after hiring Member 14 (within four weeks). New rows and columns introduced by Member 14 are highlighted in blue, and cells whose values changed relative to the pre-hiring MDM are highlighted in red.}
\label{fig:mdm-after}
\end{subfigure}
\vfill
\begin{subfigure}{\linewidth}
\includegraphics[page=8, width=0.99\columnwidth]{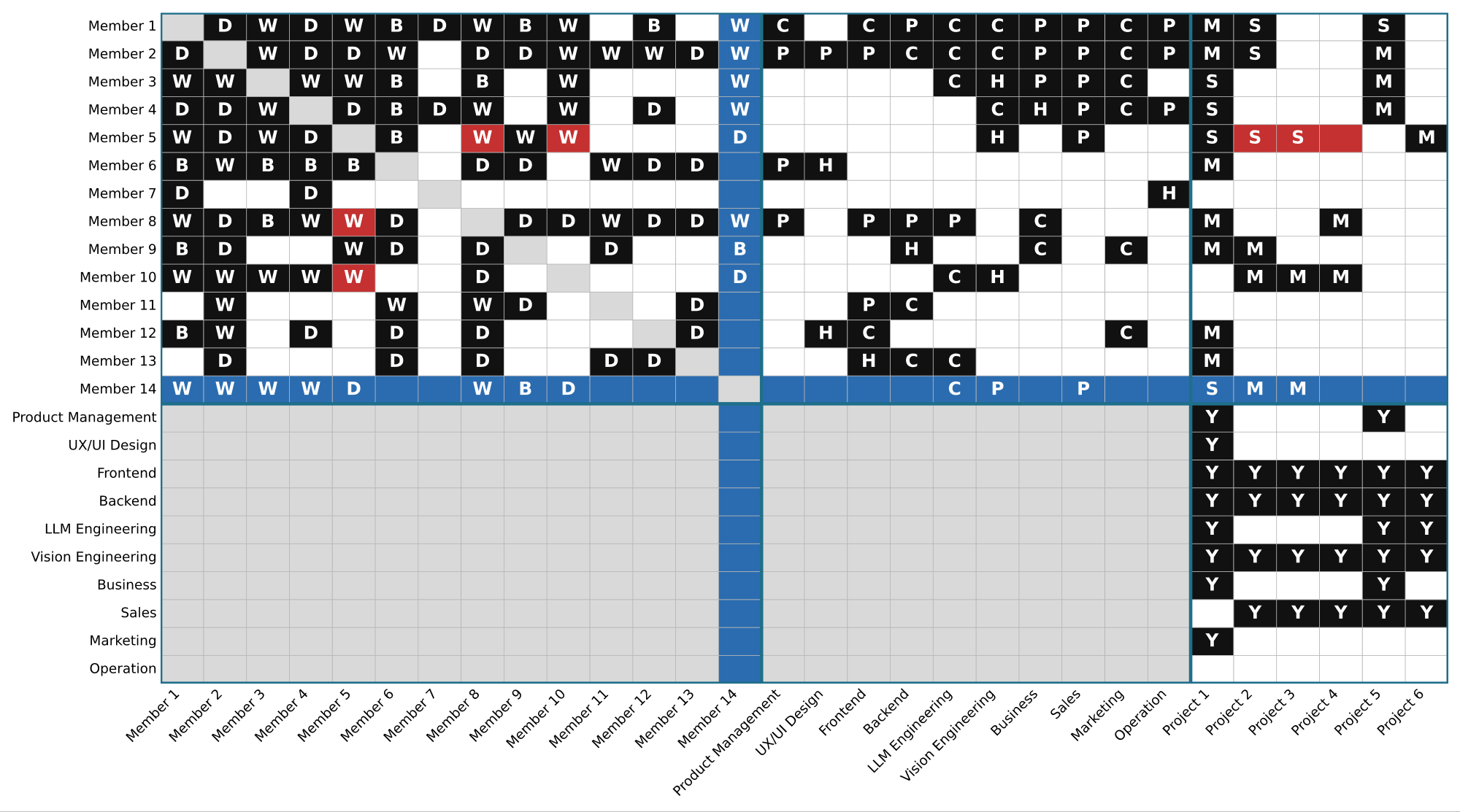}
\caption{Projected MDM three months after hiring, reflecting Member 14's full integration and the corresponding redistribution of Member 5's communication and project load.}
\label{fig:mdm-after-3month}
\end{subfigure}
\caption{MDM updates following the hiring of Member 14.}
\label{fig:mdm-updates}
\end{figure*}

\subsubsection{MDM Updates and Qualitative Changes}
\label{subsubsec:qualitative-delta}

The two updated MDMs are shown in Figure~\ref{fig:mdm-updates}. Figure~\ref{fig:mdm-after} captures the state immediately after hiring. As described in the scenario, Member 5 communicates with Member 14 daily during the handover period, and although Member 14 has assumed the primary role on P3, Member 5 remains on the project as a sub-role contributor. The remaining communication and project assignments are otherwise inherited from the pre-hiring MDM.

Figure~\ref{fig:mdm-after-3month} shows the projected state three months later. As Member 14 fully adapts to the team's workflow, Member 5 is removed from P3 entirely. Member 14 also takes over the primary role on Project~2 (P2), with Member 5 transitioning to a sub-role on that project as the handover proceeds. As workload redistribution stabilizes, Member 5's previous sub-role on Project~4 (P4) is also absorbed by the project's two main contributors, removing Member 5 from P4 altogether. Communication patterns adjust accordingly: the daily exchanges between Member 5 and Members 8 and 10---driven primarily by their shared participation in multiple projects---decrease in frequency, and the corresponding edges in the Members--Members interface are updated from Daily to Weekly.

\subsubsection{Quantitative Changes}
\label{subsubsec:quantitative-delta}

The base metrics and composite indices were recomputed at both time points using the same procedure as in Section~\ref{subsubsec:quantitative-results}. Tables~\ref{table:scores-imm} and~\ref{table:scores-3m} report the normalized base metrics and composite indices for all 14 members at the immediate post-hiring snapshot and the three-month projection, respectively, in the same format as Table~\ref{table:scores}.

\begin{table*}[h!]
\centering
\caption{Normalized base metrics and composite indices for all 14 team members after hiring Member 14. Original values are shown in parentheses for the basic metrics, whereas rankings are shown in parentheses for the composite indices. Target members highlighted.}
\label{table:scores-updated}
\begin{subtable}[t]{0.49\textwidth}
\centering
\caption{Immediately after hiring.}
\label{table:scores-imm}
\rowcolors{2}{white}{gray!10}
\scriptsize
\setlength{\tabcolsep}{3pt}
\begin{tabular}{l|cccc|cc}
\rowcolor{gray!25}
\textbf{Mem.} & \textbf{Comm} & \textbf{SL} & \textbf{Proj} & \textbf{Gap} & \textbf{WI} & \textbf{VI} \\
1  & 0.67 (22) & 0.89 (57) & 0.50 (6)  & 0.59 (1.44) & 0.56 (7)  & 0.73 (4) \\
2  & 1.00 (30) & 1.00 (64) & 0.75 (9)  & 0.36 (0.89) & 0.79 (3)  & 0.93 (1) \\
3  & 0.33 (14) & 0.86 (55) & 0.42 (5)  & 0.73 (1.78) & 0.42 (10) & 0.62 (6) \\
4  & 0.75 (24) & 0.50 (32) & 0.42 (5)  & 0.86 (2.11) & 0.56 (6)  & 0.53 (8) \\
\textbf{5}  & \textbf{0.67 (22)} & \textbf{0.70 (45)} & \textbf{0.92 (11)} & \textbf{0.95 (2.33)} & \textbf{0.85 (1)}  & \textbf{0.76 (3)} \\
6  & 0.58 (20) & 0.17 (11) & 0.33 (4)  & 0.92 (2.25) & 0.47 (9)  & 0.30 (11) \\
7  & 0.00 (6)  & 0.00 (0)  & 0.00 (0)  & 0.00 (0.00) & 0.00 (14) & 0.00 (14) \\
8  & 1.00 (30) & 0.83 (53) & 0.67 (8)  & 0.64 (1.56) & 0.76 (4)  & 0.81 (2) \\
9  & 0.42 (16) & 0.52 (33) & 0.67 (8)  & 1.00 (2.44) & 0.62 (5)  & 0.54 (7) \\
10 & 0.46 (17) & 0.52 (33) & 1.00 (12) & 0.92 (2.25) & 0.83 (2)  & 0.65 (5) \\
11 & 0.25 (12) & 0.38 (24) & 0.00 (0)  & 0.00 (0.00) & 0.07 (13) & 0.24 (13) \\
12 & 0.38 (15) & 0.19 (12) & 0.33 (4)  & 0.97 (2.38) & 0.41 (11) & 0.27 (12) \\
13 & 0.38 (15) & 0.61 (39) & 0.33 (4)  & 0.97 (2.38) & 0.41 (11) & 0.48 (10) \\
\textbf{14} & \textbf{0.46 (17)} & \textbf{0.56 (36)} & \textbf{0.42 (5)}  & \textbf{0.91 (2.22)} & \textbf{0.48 (8)}  & \textbf{0.50 (9)} \\
\hline
\textbf{Mean} & 0.55 & 0.55 & 0.48 & 0.66 & 0.52 & 0.53 \\
\end{tabular}
\end{subtable}
\hfill
\begin{subtable}[t]{0.49\textwidth}
\centering
\caption{Three months after hiring.}
\label{table:scores-3m}
\rowcolors{2}{white}{gray!10}
\scriptsize
\setlength{\tabcolsep}{3pt}
\begin{tabular}{l|cccc|cc}
\rowcolor{gray!25}
\textbf{Mem.} & \textbf{Comm} & \textbf{SL} & \textbf{Proj} & \textbf{Gap} & \textbf{WI} & \textbf{VI} \\
1  & 0.67 (22) & 0.89 (57) & 0.50 (6)  & 0.59 (1.44) & 0.56 (8)  & 0.73 (3) \\
2  & 1.00 (30) & 1.00 (64) & 0.75 (9)  & 0.36 (0.89) & 0.79 (2)  & 0.93 (1) \\
3  & 0.33 (14) & 0.86 (55) & 0.42 (5)  & 0.73 (1.78) & 0.42 (10) & 0.62 (6) \\
4  & 0.75 (24) & 0.50 (32) & 0.42 (5)  & 0.86 (2.11) & 0.56 (7)  & 0.53 (9) \\
\textbf{5}  & \textbf{0.58 (20)} & \textbf{0.70 (45)} & \textbf{0.58 (7)}  & \textbf{0.95 (2.33)} & \textbf{0.62 (6)}  & \textbf{0.64 (4)} \\
6  & 0.58 (20) & 0.17 (11) & 0.33 (4)  & 0.92 (2.25) & 0.47 (9)  & 0.30 (11) \\
7  & 0.00 (6)  & 0.00 (0)  & 0.00 (0)  & 0.00 (0.00) & 0.00 (14) & 0.00 (14) \\
8  & 0.96 (29) & 0.83 (53) & 0.67 (8)  & 0.64 (1.56) & 0.75 (3)  & 0.81 (2) \\
9  & 0.42 (16) & 0.52 (33) & 0.67 (8)  & 1.00 (2.44) & 0.62 (5)  & 0.54 (8) \\
10 & 0.42 (16) & 0.52 (33) & 1.00 (12) & 0.92 (2.25) & 0.82 (1)  & 0.64 (5) \\
11 & 0.25 (12) & 0.38 (24) & 0.00 (0)  & 0.00 (0.00) & 0.07 (13) & 0.24 (13) \\
12 & 0.38 (15) & 0.19 (12) & 0.33 (4)  & 0.97 (2.38) & 0.41 (11) & 0.27 (12) \\
13 & 0.38 (15) & 0.61 (39) & 0.33 (4)  & 0.97 (2.38) & 0.41 (11) & 0.48 (10) \\
\textbf{14} & \textbf{0.46 (17)} & \textbf{0.56 (36)} & \textbf{0.75 (9)} & \textbf{0.91 (2.22)} & \textbf{0.68 (4)} & \textbf{0.60 (7)} \\
\hline
\textbf{Mean} & 0.51 & 0.55 & 0.48 & 0.66 & 0.51 & 0.52 \\
\end{tabular}
\end{subtable}
\end{table*}

In the immediate post-hiring snapshot (Table~\ref{table:scores-imm}), Member 14 enters the organization with WI~$=0.48$ and VI~$=0.50$. Member 5's Workload Index decreases marginally from 0.87 to 0.85, as the immediate sub-role transition on P3 already outweighs the additional daily communication required to mentor Member 14, while Member 5's Value Index changes only slightly over the same period. The transient cost of onboarding therefore does not push Member 5's burden above the pre-hiring level, even before Member 14 is fully productive.

\begin{figure*}[h!]
\centering
\includegraphics[page=9, width=0.98\textwidth]{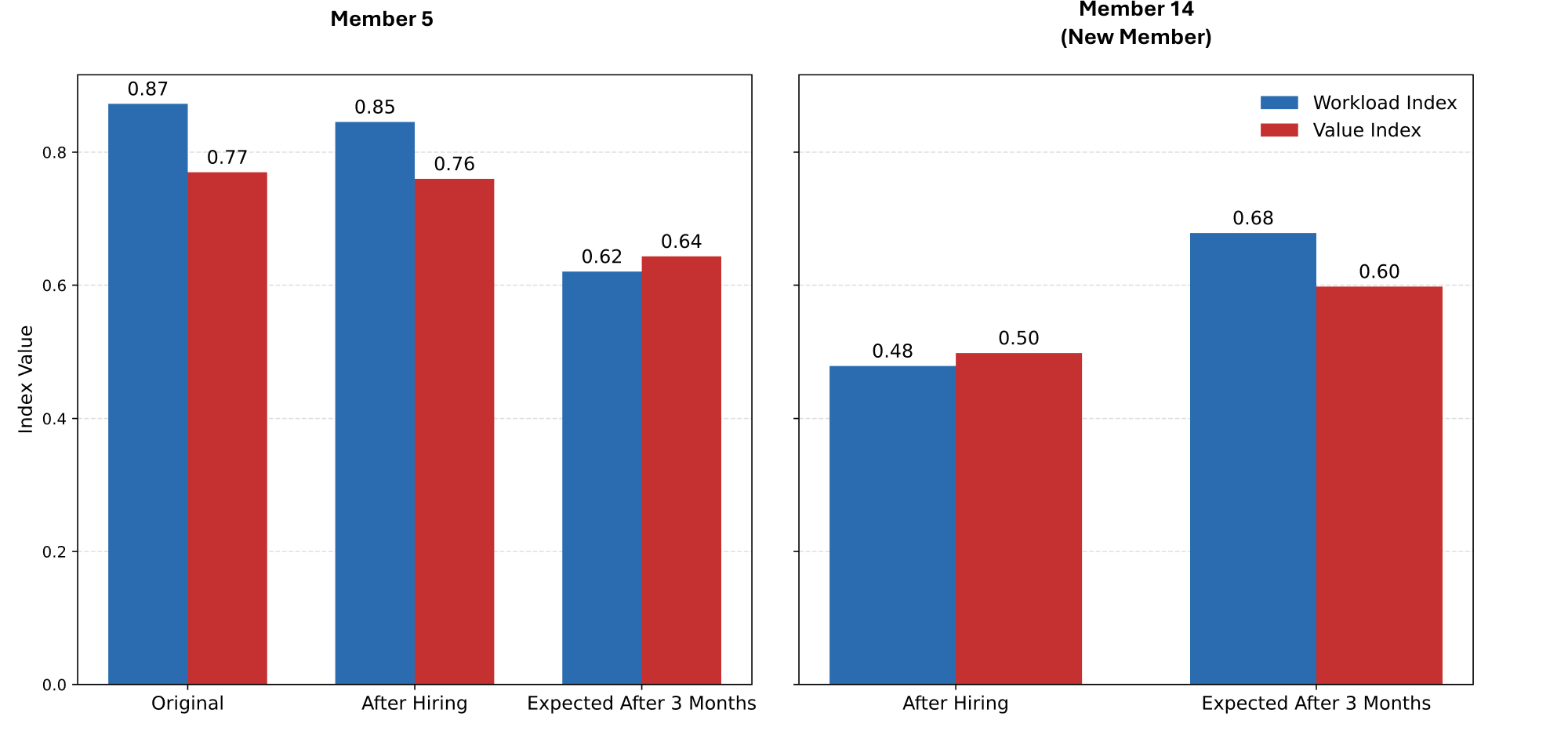}
\caption{Longitudinal Workload and Value Index changes for Members 5 and 14 across the three reference points.}
\label{fig:index-changes}
\end{figure*}

By the three-month projection (Table~\ref{table:scores-3m}), Member 14's expanded role---taking over the primary lead on P2 in addition to P3---raises Member 14's WI to 0.68 and VI to 0.60. Correspondingly, Member 5's communication and project-specific burden are substantially alleviated. As illustrated in Figure~\ref{fig:index-changes}, Member 5's Workload Index drops sharply once Member 14 reaches full operational capacity, meeting the original hiring objective.

\begin{figure*}[h!]
\centering
\begin{subfigure}{0.32\textwidth}
    \centering
    \includegraphics[page=10, width=\linewidth]{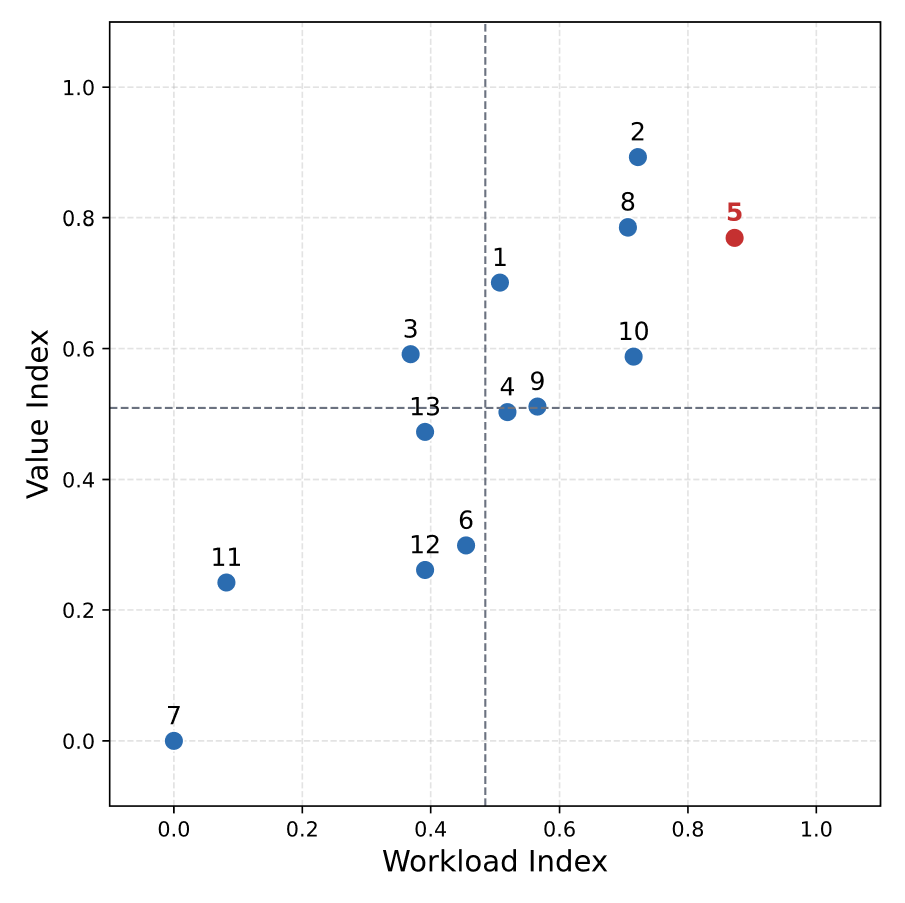}
    \caption{Pre-hiring}
\end{subfigure}
\hfill
\begin{subfigure}{0.32\textwidth}
    \centering
    \includegraphics[page=11, width=\linewidth]{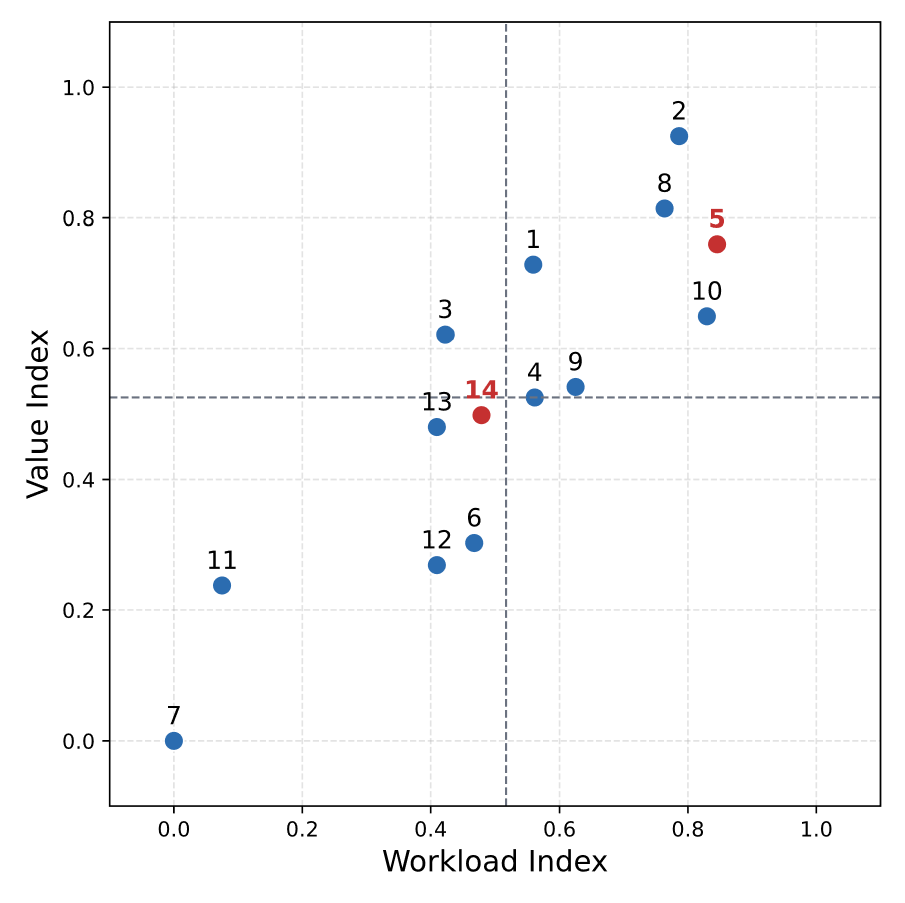}
    \caption{Immediately after hiring}
\end{subfigure}
\hfill
\begin{subfigure}{0.32\textwidth}
    \centering
    \includegraphics[page=12, width=\linewidth]{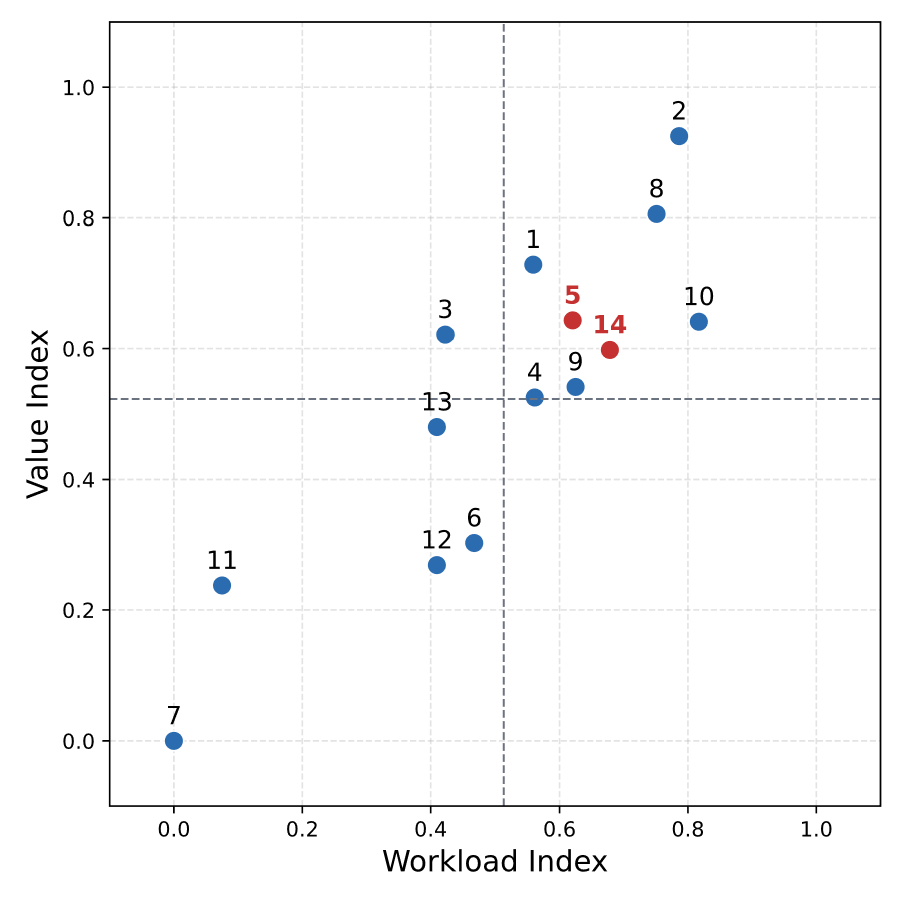}
    \caption{Three months after hiring}
\end{subfigure}
\caption{Comparative quadrant analysis at the three reference points. Members 5 and 14, the focal members of the hiring scenario, are highlighted in red; all other members are shown in blue.}
\label{fig:quadrant-after}
\end{figure*}

Figure~\ref{fig:quadrant-after} compares the three quadrant snapshots. Member 5's Workload Index decreases from 0.87 to 0.62 over the three-month horizon, an overall reduction of approximately 29\%, representing a substantial risk reduction within the Q1 region. Specifically, Member 5's Proj raw score is halved (from 14.0 to 7.0) as the member transitions from leading three primary projects to one while maintaining supporting roles on the remaining engagements. The Comm raw score also returns to a near-baseline level (from 22.0 in the immediate snapshot back to 20.0) as the daily handover communication transitions to weekly cadence after Member 14's adaptation. As a result, Member 5 remains a high-value contributor in Q1 but operates at a significantly more sustainable level, freeing cognitive bandwidth for further value-adding activities without the previous burnout risk.

Member 14's own trajectory likewise indicates successful integration. Entering at WI~$=0.48$ and VI~$=0.50$ immediately after hiring, Member 14 reaches WI~$=0.68$ (rank 4 of 14) and VI~$=0.60$ (rank 7 of 14) at the three-month projection, settling within Q1 as a core contributor rather than a peripheral hire.

More fundamentally, the intervention directly resolves Member 5's structural concentration. The number of concurrent primary roles drops from three to one over the three-month horizon, with P3 transferred to Member 14 immediately and P2 transferred during the subsequent stabilization. At the same time, the framework surfaces Member 10---who continues to hold three primary roles (P2, P3, P4) and rises to the highest Workload Index in the three-month snapshot---as the natural next candidate for intervention, illustrating how the same MDM-based pipeline can be applied iteratively as conditions evolve.

\subsubsection{Results Summary and Discussion}
\label{subsubsec:post-hiring-discussion}
Beyond these quantitative results, the case study yielded several practical insights for startup practitioners. First, multi-dimensional understanding reveals hidden dynamics: Member 5's high workload arose from a complex intersection of overlapping projects requiring specialized skills, high communication overhead, and skill gaps in adjacent technical areas---complexity that single-dimensional analysis would have missed. Second, quantitative metrics accelerate decision-making; presenting quantitative data during the hiring approval process reduced decision-making time from 2--3 weeks to approximately 1 week, with executives finding the quadrant visualization particularly effective. Third, framework implementation is feasible for small teams, as initial MDM construction required approximately 10 hours while subsequent maintenance required only 1--2 hours of monthly updates. Fourth, proactive intervention prevents costly disruptions: by identifying Member 5's burnout risk before performance degradation or departure, the organization avoided potentially severe consequences. Finally, MDM analysis complements but does not replace judgment, as the framework provides structured data to inform decisions but cannot capture qualitative factors such as team culture and organizational fit.

\section{Conclusions}
\label{sec:conclusions}

This study proposed an MDM-based HR decision support framework for small firms and startups. The framework derives four base metrics---Comm, Proj, SL, and Gap---from multi-domain matrix data and synthesizes them into a WI and a VI. A quadrant analysis based on these two composite indices enables decision-makers to classify team members by workload--value profiles and identify those requiring intervention. The framework is designed to be generalizable: domain definitions, interface measurements, and weight parameters are adaptable to different organizational contexts, while the core metric formulations and quadrant interpretation remain consistent.

The framework was demonstrated through a case study at Planby Technologies, Inc., a 13-member early-stage technology startup. The analysis identified a key team member facing unsustainable workload ($\mathrm{WI}=0.87$), which informed a targeted hiring decision. The subsequent MDM reconstruction revealed a 29\% reduction in that member's WI and resolved the concurrent triple primary-role concentration, while surfacing the next priority candidate for intervention, demonstrating the framework's ability to support iterative HR decision-making.

Several limitations and corresponding directions for future work merit mention. First, the current implementation relies on manual data collection and analysis, which may limit scalability as organizations grow or update cycles shorten. AI-based automation offers a promising direction for addressing this issue, both for initial MDM construction---by extracting member-skill-project relationships from organizational artifacts such as HR records, communication logs, and project management tools---and for ongoing updates that support near-real-time monitoring as conditions evolve. In addition, the current case study focuses on a single organization and primarily examines a hiring decision. Future work could extend the proposed framework to diverse organizational contexts and other HR decision domains, including promotion, compensation, training prioritization, and team restructuring, thereby evaluating its generalizability and broader applicability to small-organization HR practice. Lastly, this study focuses primarily on member-level metrics and HR analytics for Members ($\mathcal{N}_1$). Future research could shift or extend the analytical focus to Projects ($\mathcal{N}_3$), developing project-level metrics that evaluate organizational project portfolios, identify structurally vulnerable or resource-intensive projects, and derive corresponding managerial actions.

\section*{Acknowledgements}
\label{sec:acknowledgements}

The authors would like to thank the Planby Technologies, Inc. team members for their participation in data collection and interviews.

\section*{Conflict of Interest}
\label{sec:conflict}

Author D. Lee is affiliated with Planby Technologies, Inc., which served as the case study organization. Data collection and analysis procedures were reviewed by the co-authors to ensure objectivity, and organizational consent for publication was obtained. All other authors declare no conflict of interest.

\selectlanguage{english}
\bibliography{references.bib}

\end{document}